\documentclass[12pt,reqno]{amsart}

\makeatletter

\usepackage{amsmath,amssymb,comment,doublespace,euscript}

\theoremstyle{plain}    
\newtheorem{thm}{Theorem}[section]
\numberwithin{figure}{section} %% Comment out for sequentially-numbered
\theoremstyle{plain}    
\newtheorem{cor}[thm]{Corollary} %%Delete [thm] to re-start numbering
\newtheorem{lemma}[thm]{Lemma} %%Delete [thm] to re-start numbering
\newtheorem{lemmadef}[thm]{Lemma and Definition} %%Delete [thm] to re-start numbering
\newtheorem{prop}[thm]{Proposition}
\newtheorem{defi}[thm]{Definition}
\newtheorem{notation}[thm]{Notation}
\newtheorem{example}[thm]{Example}
\theoremstyle{remark}    
\newtheorem{rem}[thm]{Remark}

\newtheorem{claim}{Claim}[thm]

\makeatother

\newcount\theTime
\newcount\theHour
\newcount\theMinute
\newcount\theMinuteTens
\newcount\theScratch
\theTime=\number\time
\theHour=\theTime
\divide\theHour by 60
\theScratch=\theHour
\multiply\theScratch by 60
\theMinute=\theTime
\advance\theMinute by -\theScratch
\theMinuteTens=\theMinute
\divide\theMinuteTens by 10
\theScratch=\theMinuteTens
\multiply\theScratch by 10
\advance\theMinute by -\theScratch

\def\today{{\number\day\space
 \ifcase\month\or
  January\or February\or March\or April\or May\or June\or
  July\or August\or September\or October\or November\or December\fi
 \space\number\year}}

\newcommand\btimes{\displaystyle\operatornamewithlimits\times}

\newcommand\Cpx{{\mathbf C}}

\newcommand\Db{{\mathbf D}}

\newcommand\diag{\text{\rm diag}}

\newcommand\dif{\mbox{\it d}}

\newcommand\Dwt{{\widetilde{D}}}

\newcommand\dt{{\tilde{d}}}

\newcommand\Eb{{\mathbf E}}

\newcommand\eqdef{{\;\overset{\mbox{\scriptsize def}}{=}}}

\newcommand\GRM{{\text{\rm GRM}}}

\newcommand\HEu{{\EuScript H}}                   % requires package euscript

\newcommand\Hil{{\EuScript H}}                   % requires package euscript

\newcommand\HURM{{\operatorname{HURM}}}

\newcommand\id{\mathrm{id}}

\newcommand\ImagPart{{\mathrm{Im}\;}}

\newcommand\LEu{{\EuScript L}}                   % requires package euscript

\newcommand\lspan{\mathrm{span}\,}

\newcommand\Mc{{\mathcal{M}}}

\newcommand\MEu{{\EuScript M}}                   % requires package euscript

\newcommand\mut{{\tilde{\mu}}}

\newcommand\Nats{{\mathbf N}}

\newcommand\Natsz{{\Nats\cup\{0\}}}

\newcommand\nm[1]{\|#1\|}

\newcommand\nut{{\tilde{\nu}}}

\newcommand\otdt{\otimes\cdots\otimes}

\newcommand\Pt{{\widetilde P}}

\newcommand\Rc{{\mathcal{R}}}

\newcommand\RealPart{{\mathrm{Re}\;}}

\newcommand\Reals{{\mathbf R}}

\newcommand\restrict{{\upharpoonright}}

\newcommand\SGRM{{\text{\rm SGRM}}}

\newcommand\sign{{\text{\rm sign}}}

\newcommand\St{{\widetilde S}}

\newcommand\Tcirc{{\mathbf T}}

\newcommand\TEu{{\EuScript T}}

\newcommand\tr{{\mathrm{tr}}}

\newcommand\Tr{{\mathrm{Tr}}}

\newcommand\UEu{{\EuScript U}}                   % requires package euscript

\newcommand\UTGRM{{\text{\rm UTGRM}}}

\newcommand\Vt{{\widetilde V}}

\newcommand\Wt{{\widetilde W}}

\newcommand\zbar{{\overline z}}

\newcommand\Zwt{{\widetilde{Z}}}

\begin{document}
 
\pagestyle{headings}

\begin{spacing}{1.2}

\title{Invariant Subspaces of Voiculescu's Circular Operator}
 
\author{Ken Dykema and Uffe Haagerup}

%\date{\today}

\maketitle
 
%\hfil \timeanddate \hfil
 
%\markboth{\tiny Invariant, \timeanddate}{\tiny Invariant, \timeanddate}

\section{Introduction}

The invariant subspace problem relative to a von Neumann algebra $\Mc\subseteq B(\HEu)$ asks whether
every operator $T\in\Mc$ has a proper, nontrivial invariant subspace $\HEu_0\subseteq\HEu$ such that
the orthogonal projection $p$ onto $\HEu_0$ is an element of $\Mc$;
equivalently, it asks whether there is a projection $p\in\Mc$, $p\notin\{0,1\}$, such that $Tp=pTp$.
Even when $\Mc$ is a II$_1$--factor, this invariant subspace problem remains open.
In this paper we show that the circular operator and each circular free Poisson operator (defined below) has a continuous
family of invariant subspaces relative to the von Neumann algebra it generates.
These operators arise naturally in free probability theory, (see the book~\cite{VDN}), and each generates the
von Neumann algebra II$_1$--factor
$L(F_2)$ associated to the nonabelian free group on two generators.

Given a von Neumann algebra $\Mc$ with normal faithful state $\phi$, a circular operator is $y=(x_1+ix_2)/\sqrt2\in\Mc$, where
$x_1$ and $x_2$ are centered semicircular elements having the same second moments and that are free with respect to $\phi$.
For specificity, we will always take cicular elements to have the normalization $\phi(y^*y)=1$,
which is equivalent to $\phi(x_i^2)=1$.
Voiculescu found~\cite{Voiculescu:RandMat} a matrix model for a circular element, showing that if $X(n)$ is a random matrix
whose entries are i.i.d.\ complex $(0,1/n)$--Gaussian random variables then $X(n)$ converges in $*$--moments
as $n\to\infty$ to a circular element, meaning that
\begin{equation*}
\lim_{n\to\infty}\tau_n(X(n)^{\epsilon(1)}X(n)^{\epsilon(2)}\cdots X(n)^{\epsilon(k)})=\phi(y^{\epsilon(1)}y^{\epsilon(2)}\cdots y^{\epsilon(k)})
\end{equation*}
for every $k\in\Nats$ and for every choice of
$\epsilon(j)$ being ``$*$'' or no symbol, where $\tau_n$ is the expectation of the normalized trace and where
$y$ is a circular element.
Using the matrix model, Voiculescu showed that if $(y_{ij})_{1\le i,j\le N}$ is a $*$--free family of circular elements in a von Neumann algebra $\Mc$
with respect to a normal faithful state $\phi$, then the matrix $y=\frac1{\sqrt N}(y_{ij})_{1\le i,j\le N}\in M_N(\Mc)$ is circular with respect to
the state $\phi_N$ given by $\phi_N\bigl((x_{ij})_{1\le i,j\le N}\bigr)=\frac1N\sum_{i=1}^N\phi(x_{ii})$.
Furthermore, he showed that the polar decomposition of a circular operator $y$ is $y=ub$ where $u$ is a Haar unitary (i.e.\ a unitary satisfying
$\phi(u^k)=0$ for every integer $k>0$), where $b$ is a quarter--circular element,
(i.e. having moments $\phi(b^k)=\frac1\pi\int_0^2t^k\sqrt{4-t^2}$)
and where $u$ and $b$ are $*$--free.
These and results of a similar nature have been instrumental in applications of free probability to the study of the free group factors $L(F_n)$ and
related factors;
some of the first of these were~\cite{V90}, \cite{Ra92}, \cite{Dy:FreeProdR}, \cite{Ra94}, \cite{D:FreeDim}.

Voiculescu's matrix model for the circular element is the starting point for finding invariant subspaces.
Combined with a result of Dyson, it leads to upper triangular matrix models for the circular operator, namely, a sequence $Y(n)$
of upper triangluar random matrices whose $*$--moments converge to those of a circular operator.
In these models, the elements of $Y(n)$ that are above the diagonal are complex $(0,1/n)$--Gaussian random variables,
and we show that a number of different choices are possible for the diagonal entries.
From these matrix models we show for any $N\ge2$ that a circular operator can be realized as an $N\times N$ upper triangular matrix
\begin{equation}
\label{eq:utop}
\frac1{\sqrt N}\left(\begin{matrix}a_1&b_{12}&b_{13}&\cdots&b_{1,N-1}&b_{1N}\\0&a_2&b_{23}&\cdots&b_{2,N-1}&b_{2N}\\0&0&a_3&\ddots&\vdots&b_{3N}\\
\vdots&\vdots&\ddots&\ddots&\ddots&\vdots\\ 0&0&\cdots&0&a_{N-1}&b_{N-1,N}\\0&0&\cdots&0&0&a_N\end{matrix}\right)
\end{equation}
with entries in some W$^*$--noncommutative probability space, where the collection of $N(N+1)/2$ nonzero entries is $*$--free,
where the entries $b_{ij}$ lying strictly above the diagonal are circular elements
and where the entries $a_j$ on the diagonal are {\it circular free Poisson elements}, ($a_j$ being circular free Poisson of parameter $j$).
These latter are generalizations of the circular operator (in the family of R--diagonal elements introduced by Nica and Speicher~\cite{NS})
which are quite natural from the perspective of free probability theory.
\begin{defi}
Let $(A,\psi)$ be a W$^*$--noncommutative probability space with $\psi$ faithful and let $c\ge1$.
A {\em circular free Poisson} element of parameter $c$ in $(A,\psi)$ is an element of the form
$UH_c$ where
$U,\,H_c\in A$, $U$ is a Haar unitary, $H_c\ge0$, the pair $\{U,\,H_c\}$ is $*$--free
and $H_c^2$ has moments equal to those of a
free Poisson distribution\footnote{\rm 
 We must point out that the formula found on~\cite[p.\ 35]{VDN}
 for the free Poisson distribution has some errors.
 The formula for $\Rc_\mu$ found there is correct, but the formulae for $G_\mu(z)$ and for the density are incorrect.}
with parameter $c$.
Thus, letting $a=(1-\sqrt c)^2$ and $b=(1+\sqrt c)^2$,
the moments of $H_c^2$ are equal to those of the measure $\nu_c$ that is supported on $[a,b]$,
is absolutely continuous with respect to Lebesgue measure and has density
\begin{equation*}
\frac{\dif\nu_c}{\dif\lambda}(t)=\frac{\sqrt{(b-t)(t-a)}}{2\pi t}1_{[a,b]}(t). 
\end{equation*}
\end{defi}
\noindent
We hasten to point out that a circular element $z$ with normalization $\psi(z^*z)=1$
is nothing other than a circular free Poisson element of parameter $c=1$.

The spectrum of a circular free Poission element of parameter $c$ has been found by Haagerup and Larsen~\cite{HL}
to be the annulus centered at zero
with radii $\sqrt{c-1}$ and $\sqrt c$, if $c>1$, whereas the spectrum of the circular operator is the disk of radius $1$.
In the realization~\eqref{eq:utop} of the circular operator, we have that the diagonal entry $a_j$ is circular free Poisson of parameter $j$.
Therefore, the spectrum of the diagonal entry $a_j$ increases in modulus as $j$ increases, and the spectra of $a_j$ and $a_k$ overlap only
if $|j-k|\le1$.
These properties of the realization~\eqref{eq:utop} allow general techniques for upper triangular operators to be applied in order
to find invariant subspaces of the circular operator $y$.
It turns out that for every $0<r<1$ there is a unique projection $p\in\Mc=\{y\}''$ such that
\renewcommand{\labelenumi}{(\roman{enumi})}
\begin{enumerate}
\item $yp=pyp$
\item $\sigma(yp)=\{z\in\Cpx\mid|z|\le r\}$
\item $\sigma((1-p)y)=\{z\in\Cpx\mid r\le|z|\le1\}$
\end{enumerate}
where in (ii) (respectively (iii)), the spectrum is computed relative to the algebra $p\Mc p$, (respectively $(1-p)\Mc(1-p)$).

In fact, the techniques outlined above can be employed with very little extra effort
to find invariant subspaces for every circular free Poisson operator,
and the proof is presented in this generality throughout.

In~\S\ref{sec:isut}, some theory is developed proving the existence of invariant subspaces, relative to the generated von Neumann algebras,
of upper triangular operators, the spectra of whose diagonal entries satisfy certain conditions.
In~\S\ref{sec:utrandmat}, we consider upper triangular random matrices whose entries strictly above the diagonal
are i.i.d.\ complex Gaussian random variables which are independent of the diagonal entries.
The general theme of the results in~\S\ref{sec:utrandmat} is that the diagonal entries may be changed in certain
ways but that as matrix size tends to infinity, the limit $*$--moments remain the same.
In~\S\ref{sec:asympfreerandmat}, we generalize asymptotic freeness results which Voiculescu originally proved~\cite{Voiculescu:RandMat}
for Gaussian random matrices and constant diagonal matrices;
we allow the diagonal matrices to be random, subject to certain conditions.
In~\S\ref{sec:utcfP}, the random matrix results of the previous two sections together with results of Dyson and others are used
to find various upper triangular matrix models for circular free Poisson elements, and these are in turn used to find an upper
triangular realization of the same, as in~\eqref{eq:utop}.
In~\S\ref{sec:iscfP}, this upper triangular realization of the circular free Poisson element is fed into the machinery of~\S\ref{sec:isut}
to find invariant subspaces.

\section{An invariant subspace for an upper triangular operator}
\label{sec:isut}

Suppose $\Hil$ is a Hilbert space and $T:\Hil\to\Hil$ is a bounded operator.
In this section we are concerned with invariant subspaces $\Hil_0$ for $T$ such
such that the projection from $\Hil$ onto $\Hil_0$ lies in the von Neumann
algebra generated by $T$.
It is easy to see (Lemma~\ref{lemma-rinvsubsp}) that for every $r\ge0$ the set
\begin{equation}
\label{eq:HrT}
\Hil_r(T)=\overline{\{\xi\in\Hil\mid
\limsup_{k\to\infty}\nm{T^k\xi}^{1/k}\le r\}}
\end{equation}
is such an invariant subspace $\Hil_0$.
The question is, for any given operator, whether these subspaces can be other
than $\{0\}$ or $\Hil$.
We will show (Proposition~\ref{prop-uptrinv}) that the answer is yes if $T$ can
be written as an upper triangular operator,
\begin{equation*}
T=\left(\matrix*&*&*\\0&*&*\\0&0&*\endmatrix\right)
\end{equation*}
with respect to a decomposition $\Hil=\Hil_1\oplus\Hil_2\oplus\Hil_3$ under a
condition on the spectra of the elements in the upper left--hand and lower right--hand corners
of the above matrix.

\begin{lemmadef}
\label{lemma-rinvsubsp}
Let $T:\Hil\to\Hil$ be a bounded operator on a Hilbert space $\Hil$, let $r\ge0$
and let $\Hil_r=\Hil_r(T)$ be the subspace~\eqref{eq:HrT} defined above.
Then $\Hil_r$ is an invariant subspace for $T$ such that the orthogonal
projection $p_r=p_r(T)$ from $\Hil$ onto $\Hil_r$ lies in the von Neumann algebra
$\{T\}''$ generated by $T$.
\end{lemmadef}
\begin{proof}
Consider the subset
\begin{equation}
\label{eq:ErT}
E_r=E_r(T)=\{\xi\in\HEu\mid\limsup_{k\to\infty}\nm{T^k\xi}^{1/k}\le r\}
\end{equation}
which is dense in $\Hil_r$.
To see that $\Hil_r$ is a closed subspace of $\Hil$, it will suffice to show
that $E_r$ is a subspace.
Let $a\in\Cpx$ and $\xi_1,\xi_2\in E_r$.
Since $\nm{T^k(a\xi_1)}=|a|\,\nm{T^k\xi_1}$ and
$\nm{T^k(\xi_1+\xi_2)}\le2\max(\nm{T^k\xi_1},\nm{T^k\xi_2})$, it is clear that
$a\xi_1\in E_r$ and $\xi_1+\xi_2\in E_r$.
Moreover, since $\nm{T^k(T\xi)}\le\nm T\nm{T^k(\xi)}$ it is clear that $E_r$,
and hence also $\Hil_r$, is invariant for $T$.

To show that $p_r\in\{T\}''$ it will be enough to show that $Up_r=p_rU$
whenever $U$ is a unitary operator on $\Hil$ such that $UT=TU$.
Moreover, $Up_r=p_rU$ will follow once we show that $U\xi\in E_r$ for every
$\xi\in E_r$.
But this holds because
\[
\nm{T^kU\xi}=\nm{UT^k\xi}=\nm{T^k\xi}.
\]
\end{proof}

\begin{prop}
\label{prop-uptrinv}
Let $T:\Hil\to\Hil$ be a bounded operator on a Hilbert space $\Hil$.
Suppose that $e_1,e_2,e_3$ are orthogonal projections in $\Hil$ with
$e_1+e_2+e_3=1$ and that $e_1$ and $e_1+e_2$ are invariant for $T$.
This means that $T$ is upper triangular with respect to this decomposition
of $\Hil$:
\[
T=\left(\matrix
e_1Te_1&*&*\\
0&e_2Te_2&*\\
0&0&e_3Te_3\endmatrix\right).
\]
Let $r\ge0$ and suppose that
\[
\sup\{|\lambda|\mid\lambda\in\sigma(e_1Te_1)\}\le r
<\inf\{|\lambda|\mid\lambda\in\sigma(e_3Te_3)\},
\]
where $\sigma(e_jTe_j)$ denotes the spectrum of $e_jTe_j$ acting on $e_j\HEu$.
Then the invariant subspace $\Hil_r(T)$ and its projection $p_r=p_r(T)$ defined in
Lemma and Definition~\ref{lemma-rinvsubsp} satisfy
\[
e_1\le p_r\le e_1+e_2.
\]
\end{prop}
\begin{proof}
If $\xi=e_1\xi$ then $T^k\xi=(e_1Te_1)^k\xi$ so
$\nm{T^k\xi}^{1/k}\le\bigl(\nm{(e_1Te_1)^k}{\nm\xi}\bigr)^{1/k}$ while
$\lim_{k\to\infty}\linebreak\nm{(e_1Te_1)^k}^{1/k}\le r$.
This shows $e_1\le p_r$.

Suppose $\xi\in\Hil$ and $e_3\xi\ne0$.
Then $e_3T^k\xi=(e_3Te_3)^k\xi$.
As an operator on $e_3\Hil$, $e_3Te_3$ is invertible and its inverse has
spectral radius $<r^{-1}$.
Therefore
$\nm{T^k\xi}\ge\nm{(e_3Te_3)^k\xi}\ge\nm{(e_3Te_3)^{-k}}^{-1}\nm{e_3\xi}$ and
hence
$\limsup_{k\to\infty}\nm{T^k\xi}^{1/k}\ge
\lim_{k\to\infty}\nm{(e_3Te_3)^{-k}}^{-1/k}>r$,
so $\xi\not\in E_r(T)$.
This shows that $E_r(T)\perp e_3\Hil$, and therefore that $p_r\le e_1+e_2$.
\end{proof}

Invariant projections $p_r(T)$ and the dense subspaces $E_r(T)$ for some specific operators $T$ are described in~\S\ref{sec:iscfP}.

Now we show that for an element $x$ of an abstract W$^*$--algebra, the projection $p_r(x)$ is defined independently of how
the W$^*$--algebra is represented as a von Neumann algebra acting on a Hilbert space.

\begin{lemma}
\label{lem:prTot}
Let $\Hil$ and $\Hil'$ be Hilbert spaces, let $T\in B(\Hil)$ and take $r>0$.
Then
\begin{equation}
\label{eq:prTot}
p_r(T)\otimes1_{\Hil'}=p_r(T\otimes1_{\Hil'}).
\end{equation}
\end{lemma}
\begin{proof}
Let $E_r(T)$ be given by~\eqref{eq:ErT} and let
\begin{equation*}
E_r(T\otimes1)=\{w\in\Hil\otimes\Hil'\mid\limsup_{k\to\infty}\nm{(T^k\otimes1_{\Hil'})w}^{1/k}\le r\}.
\end{equation*}
Letting $\odot$ denote the algebraic tensor product of vector spaces, we clearly have
$E_r(T)\odot\Hil'\subseteq E_r(T\otimes1)$, so the inequality $\le$ holds in~\eqref{eq:prTot}.
Given a unit vector $\eta\in\Hil'$ let $V_\eta:\Hil\to\Hil\otimes\Hil'$ be $V_\eta(\xi)=\xi\otimes\eta$.
Then $TV_\eta^*=V_\eta^*(T\otimes1_{\Hil'})$, so if $w\in E_r(T\otimes1)$ then $V_\eta^*w\in E_r(T)$ for every $\eta$.
Therefore $E_r(T\otimes 1)\subseteq\Hil_r(T)\otimes\Hil'$ and $\ge$ holds in~\eqref{eq:prTot}.
\end{proof}

\begin{lemmadef}
\label{lemmadef:invvN}
If $\Mc$ is a von Neumann algebra, if $\HEu$ is a Hilbert space and if $\pi:\Mc\to B(\HEu)$ is a normal, faithful $*$--representation
then given $x\in\Mc$ we have, for $r\ge0$, the projection
\begin{equation*}
p_r(\pi(x))\in\{\pi(x)\}''\subseteq\pi(\Mc).
\end{equation*}
Let us denote by $p_r(x)\in\Mc$ the element so that $\pi(p_r(x))=p_r(\pi(x))$.
Then $p_r(x)$ is independent of the choice of $\HEu$ and $\pi$.
\end{lemmadef}
\begin{proof}
Let $\HEu'$ be a Hilbert space and $\pi':\Mc\to B(\HEu')$ a normal, faithful $*$--representation.
Using~\cite[Theorem 1.4.3]{Dix}, one finds a Hilbert space $\HEu''$ such that the representations $\pi\otimes1_{\HEu''}$
and $\pi'\otimes1_{\HEu''}$ are unitarily equivalent, via a unitary $U:\HEu\otimes\HEu''\to\HEu'\otimes\HEu''$.
Now applying Lemma~\ref{lem:prTot} twice, we have
\begin{align*}
U^*\bigl(\pi'(p_r(x))\otimes1\bigr)U&=\pi(p_r(x))\otimes1=p_r(\pi(x))\otimes1=p_r(\pi(x)\otimes1)= \\
&=p_r\bigl(U^*(\pi'(x)\otimes1)U\bigr)=U^*p_r\bigl(\pi'(x)\otimes1\bigr)U=U^*\bigl(p_r(\pi'(x))\otimes1\bigr)U.
\end{align*}
Hence $\pi'(p_r(x))=p_r(\pi'(x))$.
\end{proof}

\section{Upper triangular random matrices}
\label{sec:utrandmat}

In this section, we consider upper triangular random matrices whose entries strictly above the diagonal are i.i.d.\ Gaussian,
and we prove that the diagonal entries can be modified in various ways without affecting the limiting $*$--moments as matrix size increases without bound.
Throughout this paper, we consider random matrices whose entries have moments of all orders.
Thus, let $(\Omega,\mu)$ be a usual probability space, let
$\LEu=\bigcap_{1\le p<\infty}L^p(\mu)$ and consider the expectation
$E(f)=\int f d\mu$.
If $S_1,S_2\subseteq\LEu$ are sets of random variables, we say that $S_1$ and
$S_2$ are {\it independent sets} if the two $\sigma$--algebras generated by
$\{f^{-1}(A)\mid f\in S_i,\,A\mbox{ \ Borel subset of }\Cpx\}$ ($i=1,2$) are
independent with respect to $\mu$, and similarly for families of sets of random variables.
The $*$--algebra of $n\times n$ random matrices is $\MEu_n=\LEu\otimes M_n(\Cpx)$
and has the trace $\tau_n=E\otimes\tr_n$, where $\tr_n$ is the trace on
$M_n(\Cpx)$ normalized so that $\tr_n(1)=1$.
We fix a system of matrix units in $M_n(\Cpx)$, denoted by
$(e(i,j;n))_{1\le i,j\le n}$.

\begin{notation}
\label{notation:UTGRM}
Let $\sigma^2>0$ and $n\in\Nats$.
\renewcommand{\labelenumi}{(\roman{enumi})}
\begin{enumerate}
\item On $\Cpx$, by {\em Lebesgue measure} we shall mean $\dif(\RealPart z)\dif(\ImagPart z)$, i.e. normalized so that the unit disk has measure $\pi$.
On the space $M_n(\Cpx)$ of $n\times n$ complex matrices, Lebesgue measure shall mean the product of Lebesgue measure on each of the $n^2$ complex entries.
On the space $M_n^{s.a.}$ of self--adjoint complex $n\times n$ matrices, Lebesgue measure shall mean
the product of Lebesgue measure on each of the $n(n-1)/2$ complex
entries strictly above the diagonal and Lebesgue measure on each of the $n$ real diagonal entries.

\item We say that a random variable $a\in\LEu$ is a {\em complex $(0,\sigma^2)$--Gaussian}
if $\RealPart a$ and $\ImagPart a$ are independent real Gaussian random variables each having first moment $0$ and second moment $\sigma^2/2$.
Thus $E(a)=0$, $E(|a|^2)=\sigma^2$ and $a$ has density $(\pi\sigma^2)^{-1}e^{-|z|^2/\sigma^2}$ with respect to Lebesgue measure.

\item
Given $T\in\MEu_n$, we will write $T\in\UTGRM(n,\sigma^2)$, (the acronym is for``upper triangular Gaussian random matrix'')
if the entries $t_{ij}$ of $T$ ($1\le i,j\le n$) satisfy that $t_{ij}=0$ whenever $i\ge j$ and that $(t_{ij})_{1\le i<j\le n}$ is an independent family
of random variables, each of which is complex $(0,\sigma^2)$--Gaussian.
\end{enumerate}
\end{notation}

Our first result is that if $T(n)\in\UTGRM(n,1/n)$, if $D(n)\in\MEu_n$ is a diagonal random matrix
that is independent from $T(n)$ and if the joint distribution of the diagonal entries of $D(n)$
is permutation invariant then in the limit as $n\to\infty$, the $*$--moments  of $T(n)+D(n)$
depend only on the marginal $*$--distributions of finite sets of the diagonal entries of $D(n)$.

\begin{thm}
\label{thm:utgrm}
For every $n\in\Nats$ let $T(n)\in\UTGRM(n,\frac1n)$ and let
\begin{align*}
D_1(n)&=\sum_{i=1}^na(i;n)\otimes e(i,i;n)\in\MEu_n \\
D_2(n)&=\sum_{i=1}^nb(i;n)\otimes e(i,i;n)\in\MEu_n
\end{align*}
be diagonal random matrices such that
$T(n)$ and $D_1(n)$ are independent matrix--valued random variables
and $T(n)$ and $D_2(n)$ are independent matrix--valued random variables.
Let $d\nu_\iota(\lambda_1,\ldots,\lambda_n)$ be the joint distribution of the
diagonal entries of $D_\iota$.
Assume that $d\nu_\iota$ is invariant under all permutations of its $n$ variables,
($\iota=1,2$).
Suppose that the marginal $*$--distributions of the diagonal entries of $D_1(n)$ are arbitrarily close to those of $D_2(n)$ as $n\to\infty$;
namely suppose that
\pagebreak[0]
\begin{equation}
\label{eq:abmomhyp}
\begin{split}
\forall p\in\Nats\quad\forall r_1,s_1,\ldots,r_p,&s_p\in\Natsz, \\
\lim_{n\to\infty}\bigg(&
E(a(1;n)^{r_1}\overline{a(1;n)}^{s_1}\cdots
a(p;n)^{r_p}\overline{a(p;n)}^{s_p}) \\
&-E(b(1;n)^{r_1}\overline{b(1;n)}^{s_1}\cdots
b(p;n)^{r_p}\overline{b(p;n)}^{s_p})\bigg)=0.
\end{split}
\end{equation}
Let $Z_\iota(n)=D_\iota(n)+T(n)$.
Then
\begin{align*}
\forall m\in\Nats&\quad\forall\epsilon(1),\ldots,\epsilon(m)
\in\{*,\cdot\}, \\[1.5ex]
&\lim_{n\to\infty}\big(\tau_n(Z_1(n)^{\epsilon(1)}\cdots Z_1(n)^{\epsilon(m)})
-\tau_n(Z_2(n)^{\epsilon(1)}\cdots Z_2(n)^{\epsilon(m)})\big)=0,
\end{align*}
where we take $\epsilon(j)=\cdot$ to mean $\epsilon(j)$ is ``not $*$'' or ``no symbol.''
Therefore if $Z_1(n)$ converges in $*$--moments as $n\to\infty$ then so does $Z_2(n)$ and their limit $*$--moments coincide.
\end{thm}
\begin{proof}
Write
\begin{equation*}
T(n)=\sum_{1\le i<j\le n}t(i,j;n)\otimes e(i,j;n).
\end{equation*}
Let us first fix $n\in\Nats$, $\iota\in\{1,2\}$ and let us denote $D_\iota(n)$
simply by $D$, $T(n)$ by $T$, $Z_\iota(n)$ by $Z$ and the diagonal entries of $D$ by
$d(i;n)$, ($1\le i\le n$).
Now each word in $Z$ and $Z^*$ is a sum of words in $D$, $D^*$, $T$ and $T^*$;
hence we will investigate words of the form
\begin{multline}
\label{eq:Wexpr}
W=\bigg(T^{\epsilon(1)}\cdots T^{\epsilon(l(1))}\bigg)D^{\kappa(1)}
\bigg(T^{\epsilon(l(1)+1)}\cdots T^{\epsilon(l(2))}\bigg)D^{\kappa(2)}\cdots \\
\cdots\bigg(T^{\epsilon(l(p-1)+1)}\cdots T^{\epsilon(l(p))}\bigg)D^{\kappa(p)}
\bigg(T^{\epsilon(l(p)+1)}\cdots T^{\epsilon(q)}\bigg)
\end{multline}
for arbitrary
\begin{gather*}
p,q\in\Natsz,\\
\epsilon(1),\ldots,\epsilon(q),
\kappa(1),\ldots,\kappa(p)\in\{*,\cdot\}, \\
0\le l(1)\le l(2)\le\cdots\le l(p)\le q.
\end{gather*}
Now, using the independence of T and D, we
see that
\begin{multline}
\label{eq:taunWexpr}
\tau_n(W)= \displaybreak[0] \\ 
=n^{-1}\sum_{i_1,\ldots,i_q\in\{1,\ldots,n\}}
\bigg(\prod_{j=1}^qG^{\epsilon(j)}(i_j,i_{j+1})\bigg)
E\bigg(t^{\epsilon(1)}(i_1,i_2;n)t^{\epsilon(2)}(i_2,i_3;n)\cdots
t^{\epsilon(q)}(i_q,i_{q+1};n)\bigg)\cdot \\
\cdot E\bigg(d^{\kappa(1)}(i_{l(1)};n)d^{\kappa(2)}(i_{l(2)};n)\cdots
d^{\kappa(p)}(i_{l(p)};n)\bigg),
\displaybreak[0]
\end{multline}
where we have used the convention $i_{q+1}=i_1$, where
$G^\epsilon(\cdot,\cdot)$ is defined by
\begin{align*}
G(i_j,i_{j+1})&=
\begin{cases}1&\mbox{if }i_j<i_{j+1}\\0&\mbox{otherwise}\end{cases} \displaybreak[1] \\[2ex]
G^*(i_j,i_{j+1})&=
\begin{cases}1&\mbox{if }i_j>i_{j+1}\\0&\mbox{otherwise}\end{cases}
\end{align*}
and where
\begin{align*}
t^{\epsilon(j)}(i_j,i_{j+1};n)&=
\begin{cases}t(i_j,i_{j+1};n)
&\mbox{if }\epsilon(j)=\cdot\\ [1.5ex]
\overline{t(i_{j+1},i_j;n)}&\mbox{if }\epsilon(j)=*\end{cases}\\[2.5ex]
d^{\kappa(j)}(i_j;n)&=
\begin{cases}d(i_j;n)
&\mbox{if }\kappa(j)=\cdot\\ [1.5ex]
\overline{d(i_j;n)}&\mbox{if }\kappa(j)=*\end{cases}
\end{align*}
Using that the joint distribution of $d(1;n),\ldots,d(n;n)$ is invariant under
permutation of the $n$ variables, we see that each moment
\begin{equation}
\label{eq:dmomentwild}
E\big(d^{\kappa(1)}(i_{l(1)};n)d^{\kappa(2)}(i_{l(2)};n)\cdots
d^{\kappa(p)}(i_{l(p)};n)\big)
\end{equation}
appearing in~(\ref{eq:taunWexpr}) is equal to a
moment
\begin{equation}
\label{eq:dmoment}
E\big(d(1;n)^{r(1)}\overline{d(1;n)}^{s(1)}d(2;n)^{r(2)}\overline{d(2;n)}^{s(2)}
\cdots d(p;n)^{r(p)}\overline{d(p;n)}^{s(p)}\big)
\end{equation}
where
\begin{equation}
\label{eq:rscond1}
\begin{aligned}
r(1),s(1),\ldots,r(p),s(p)&\in\Natsz,\\
r(1)+s(1)+\cdots+r(p)+s(p)&=p
\end{aligned}
\end{equation}
and by further permutation the moment~(\ref{eq:dmomentwild}) corresponds to a
unique moment of the form~(\ref{eq:dmoment}) if we make the additional
stipulation that
\begin{equation}
\label{eq:rscond2}
\begin{gathered}
r(1)+s(1)\ge r(2)+s(2)\ge\cdots\ge r(p)+s(p) \\
\mbox{and, if }r(j)+s(j)=r(j+1)+s(j+1)\mbox{ then }r(j)\ge r(j+1).
\end{gathered}
\end{equation}
Hence, rearranging the sum in~(\ref{eq:taunWexpr}) we get
\begin{multline}
\label{eq:newtaunWexpr}
\tau_n(W)
=n^{-1}\sum_{r(1),s(1),\ldots,r(p),s(p)}
E\bigg(d(1;n)^{r(1)}\overline{d(1;n)}^{s(1)}
\cdots d(p;n)^{r(p)}\overline{d(p;n)}^{s(p)}\bigg)\cdot \\
\cdot\sum_{i_1,\ldots,i_q}
\bigg(\prod_{j=1}^qG^{\epsilon(j)}(i_j,i_{j+1})\bigg)
E\bigg(t^{\epsilon(1)}(i_1,i_2;n)\cdots
t^{\epsilon(q)}(i_q,i_{q+1};n)\bigg),
\end{multline}
where the first sum is over all $r(1),s(1),\ldots,r(p),s(p)$
satisfying~(\ref{eq:rscond1}) and~(\ref{eq:rscond2}), and the second sum is
over all $i_1,\ldots,i_q\in\{1,\ldots,n\}$ such that there is a permutation,
$\sigma$, of $\{1,\ldots,n\}$ for which
\begin{equation}
\label{eq:dmom}
\begin{split}
d(\sigma(i_1);n)^{\kappa(1)}d(\sigma(i_2);n)^{\kappa(2)}&\cdots
d(\sigma(i_p);n)^{\kappa(p)}= \\
&=d(1;n)^{r(1)}\overline{d(1;n)}^{s(1)}
\cdots d(p;n)^{r(p)}\overline{d(p;n)}^{s(p)}.
\end{split}
\end{equation}

Let $W_i$ ($\iota\in\{1,2\}$) denote the word on the right--hand--side
of~(\ref{eq:Wexpr}) where $D$ is taken to be $D_\iota$.
Then
\begin{multline*}
\tau_n(W_1)-\tau_n(W_2)= \\
\begin{split}
=n^{-1}\sum_{r(1),s(1),\ldots,r(p),s(p)}
\bigg(&E\big(a(1;n)^{r(1)}\overline{a(1;n)}^{s(1)}
\cdots a(p;n)^{r(p)}\overline{a(p;n)}^{s(p)}\big) \\
&-E\big(b(1;n)^{r(1)}\overline{b(1;n)}^{s(1)}
\cdots b(p;n)^{r(p)}\overline{b(p;n)}^{s(p)}\big)\bigg)\cdot
\end{split} \\
\cdot\sum_{i_1,\ldots,i_q}
\bigg(\prod_{j=1}^qG^{\epsilon(j)}(i_j,i_{j+1})\bigg)
E\bigg(t^{\epsilon(1)}(i_1,i_2;n)\cdots
t^{\epsilon(q)}(i_q,i_{q+1};n)\bigg).
\end{multline*}
Now in order to prove the lemma it will suffice to show that
$\lim_{n\to\infty}\tau_n(W_1)-\tau_n(W_2)=0$.
Since by the hypothesis~(\ref{eq:abmomhyp}), the difference of moments of $a$
and $b$ tends to zero as $n\to\infty$, it will suffice to show that for every
$p,q\in\Nats$ and every $r(1),s(1),\ldots,r(p),s(p)\in\Natsz$, the quantity
\begin{equation}
\label{eq:tobebdd}
n^{-1}\sum_{i_1,\ldots,i_q}
\bigg(\prod_{j=1}^qG^{\epsilon(j)}(i_j,i_{j+1})\bigg)
E\bigg(t^{\epsilon(1)}(i_1,i_2;n)\cdots
t^{\epsilon(q)}(i_q,i_{q+1};n)\bigg)
\end{equation}
remains bounded as $n\to\infty$, where the sum is
over all $i_1,\ldots,i_q\in\{1,\ldots,n\}$ such that there is a permutation,
$\sigma$, of $\{1,\ldots,n\}$ for which~(\ref{eq:dmom}) holds.
But this follows from the sort of counting arguments used by Voiculescu
in~\cite{Voiculescu:RandMat}.
Indeed, for any $1\le i<i'\le n$ and for any $m,m'\in\Natsz$,
\[
E\big(t(i,i';n)^m\overline{t(i,i';n)}^{m'}\big)\ne0
\]
implies $m=m'$.
Taken together with the independence assumtion on the entries of $T$, this
shows that for any $i_1,\ldots,i_q\in\{1,\ldots,n\}$, a necessary condition so
that
\[
E\big(t^{\epsilon(1)}(i_1,i_2;n)\cdots
t^{\epsilon(q)}(i_q,i_{q+1};n)\big)\ne0
\]
is that there be a bijection, $\gamma$, from $\{1,\ldots,q\}$ to itself,
without fixed points, such that $\gamma^2=\id$ and
\begin{equation}
\label{eq:gammacond}
\forall j\in\{1,\ldots,q\},\quad
i_{\gamma(j)}=i_{j+1},\mbox{ and }i_{\gamma(j)+1}=i_j.
\end{equation}
Moreover, there is a constant, $c_1$, depending only on $q$, such that for all
$n\in\Nats$ and all choices of $i_1,\ldots,i_q$,
\begin{equation}
\label{eq:Etbnd}
\left|
E\big(t^{\epsilon(1)}(i_1,i_2;n)\cdots
t^{\epsilon(q)}(i_q,i_{q+1};n)\big)
\right|\le c_1n^{-q/2}.
\end{equation}
If $\gamma$ is a bijection of $\{1,\ldots,q\}$, let
$N(\gamma,n)$ be the number of choices of $i_1,\ldots,i_q\in\{1,\ldots,n\}$
such that~(\ref{eq:gammacond}) holds.
There are only finitely many bijections, $\gamma$, of
$\{1,\ldots,q\}$.
Hence, in light of the bound~(\ref{eq:Etbnd}), in order to show
that~(\ref{eq:tobebdd}) is bounded as $n\to\infty$, it will suffice to show
that for each bijection $\gamma$ of $\{1,\ldots,q\}$ without fixed points such
that $\gamma^2=\id$, the quantity
\begin{equation}
\label{eq:ndNgam}
n^{-1-(q/2)}N(\gamma,n)
\end{equation}
remains bounded as $n\to\infty$.
However, $N(\gamma,n)\le n^{d(\gamma)}$, where $d(\gamma)$
is the number of vertices in the quotient graph, $G_\gamma$, which is obtained
from the $q$--gon graph by identifying the $j$th and $\gamma(j)$th edges with
opposite orientations, for every $j$.
The graph $G_\gamma$ has $q/2$ edges, hence at most $1+(q/2)$ vertices,
which shows that $d(\gamma)\le1+(q/2)$ and hence that~(\ref{eq:ndNgam}) remains
bounded as $n\to\infty$.
\end{proof}

Now we work on results that let us dispense with the permutation
invariance supposed for the diagonal matrices of the previous theorem.
Let $\UEu_2$ denote the group of unitary $2\times2$ complex matrices.

\begin{lemma}
\label{lem-Uabc}
There is a Borel function, $U:\Cpx^3\to\UEu_2$ such that for all
$a,b,c\in\Cpx$,
\[
U(a,b,c)^*\begin{pmatrix}a&b\\0&c\end{pmatrix}U(a,b,c)
=\begin{pmatrix}c&b\\0&a\end{pmatrix}.
\]
\end{lemma}

\begin{proof}
If $a=c$ then let $U(a,b,c)=\begin{pmatrix}1&0\\0&1\end{pmatrix}$.
If $a\ne c$ but $b=0$ then let $U(a,b,c)=\begin{pmatrix}0&1\\1&0\end{pmatrix}$.
If $a\ne c$ and $b\ne0$ then let
\[
U(a,b,c)=\frac1{(|a-c|^2+|b|^2)^{1/2}}
\begin{pmatrix}b&\overline{(a-c)}b/\overline b\\[1ex]c-a&b\end{pmatrix}.
\]
\end{proof}

\begin{lemma}
\label{lem:permind}
Fix $n\in\Nats$, let $T\in\UTGRM(n,\frac1n)$ and let
\[
D=\sum_{i=1}^nd(i)\otimes e(i,i;n)\in\MEu_n
\]
be a diagonal random matrix.
Suppose that $T$ and $D$
are independent matrix--valued random variables.
Let $\pi$ be a permutation of $\{1,\ldots,n\}$ and let
\begin{equation}
\label{eq:Dpi}
D_\pi=\sum_{i=1}^nd(\pi(i))\otimes e(i,i;n)\in\MEu_n.
\end{equation}
Let $Z=D+T$
and $Z_\pi=D_\pi+T$.
Then $Z$ and $Z_\pi$ have the same $*$--moments with respect to $\tau_n$.
\end{lemma}

\begin{proof}
We may without loss of generality assume that $\pi$ is a transposition of
neighbors:
\begin{gather*}
\pi(k)=k+1\\
\pi(k+1)=k\\
\pi(j)=j\quad\mbox{if }j\not\in\{k,k+1\}.
\end{gather*}
We will use Lemma~\ref{lem-Uabc} to show that there is a unitary random matrix,
$V\in\MEu_n$, such that $V^*ZV$ has the same distribution as $Z_\pi$, and this
will prove the theorem.

Recall that $\Omega$ is the usual probability space underlying our random matrices $\MEu_n$.
For $\omega\in\Omega$ let $V(\omega)$ be the block diagonal matrix
\[
V(\omega)=I_{k-1}\oplus U\big(d(k)(\omega),t(k,k+1)(\omega),d(k+1)(\omega)\big)
\oplus I_{n-k-1}.
\]
By this we mean that $V(\omega)$ has $k-1$ ones down the diagonal, then a
$2\times2$ block that is the matrix $U$ from Lemma~\ref{lem-Uabc}, then
$n-k-1$ ones.
Let $x(i,j)$ denote the $(i,j)$th entry of the random matrix $V^*ZV$.
If we write
\begin{equation*}
T=\sum_{1\le i<j\le n}t(i,j)\otimes e(i,j;n)
\end{equation*}
then
\begin{alignat*}2
x(i,j)&=0&&\mbox{if }i>j \\
x(i,i)&=d(i)&&\mbox{if }i\not\in\{k,k+1\}\\
x(k,k)&=d(k+1)\\
x(k+1,k+1)&=d(k) \\
x(k,k+1)&=t(k,k+1).
\end{alignat*}
Let
\[
\begin{pmatrix}u_{11}&u_{12}\\u_{21}&u_{22}\end{pmatrix}
=U\big(d(k)(\omega),t(k,k+1)(\omega),d(k+1)(\omega)\big).
\]
If $i<k$ then
\begin{align*}
x(i,k)&=t(i,k)u_{11}+t(i,k+1)u_{21} \\
x(i,k+1)&=t(i,k)u_{12}+t(i,k+1)u_{22}
\end{align*}
and if $j>k+1$ then
\begin{align*}
x(k,j)&=\overline{u_{11}}t(k,j)+\overline{u_{21}}t(k+1,j) \\
x(k+1,j)&=\overline{u_{12}}t(k,j)+\overline{u_{22}}t(k+1,j).
\end{align*}
In order to show that $V^*ZV$ and $Z_\pi$ have the same distribution, it is
thus enough to show
\renewcommand{\labelenumi}{(\alph{enumi})}
\begin{enumerate}
\item $(x(i,j))_{1\le i<j\le n}$ is an independent family of complex $(0,1/n)$--Gaussian random variables;
\item $\{x(i,i)\mid1\le i\le n\}$ and $\{x(i,j)\mid1\le i<j\le n\}$ are independent sets of random variables.
\end{enumerate}
From the facts that
\begin{equation}
\label{eq:twosets}
\{d(i)\mid1\le i\le n\}\cup\{t(k,k+1)\}
\quad\mbox{and}\quad\{t(i,j)\mid1\le i<j\le n,\,(i,j)\ne(k,k+1)\}
\end{equation}
are independent sets of random variables, each $t(i,j)$ is complex $(0,1/n)$--Gaussian
and $U$ is everywhere unitary and is independent from the right--hand set in~\eqref{eq:twosets}, we see that
\begin{equation}
\label{eq:xijGiid}
(x(i,j))_{1\le i<j\le n,\,(i,j)\ne(k,k+1)}
\end{equation}
is an independent family of complex $(0,1/n)$--Gaussian random variables.
Moreover, the joint distribution of the family~(\ref{eq:xijGiid}) is not
changed by conditioning on the values of $d(1),\ldots,d(n),\linebreak[0]t(k,k+1)$.
Hence
\[
\{d(i)\mid1\le i\le n\}\cup\{t(k,k+1)\}
\qquad\mbox{and}\qquad\{x(i,j)\mid1\le i<j\le n,\,(i,j)\ne(k,k+1)\}
\]
are independent sets of random variables.
But this implies that~(a) and~(b) hold.
\end{proof}

\begin{lemma}
\label{lem-desymmetrize}
Let $D\in\MEu_n$ be a diagonal random matrix, let $T\in\UTGRM(n,\frac1n)$ and let $Z=D+T$.
Let $\mu$ be the joint distribution of the $n$ random variables,
$d(1),d(2),\ldots,d(n)$, in that order.
For every permutation
$\pi$ of $\{1,\ldots,n\}$ let $\mu_\pi$ be the joint distribution of the random
variables, $d(\pi(1)),d(\pi(2)),\ldots,d(\pi(n))$, in that order.
Let $A$ be a nonempty set of permutations of $\{1,\ldots,n\}$ and let $|A|$ denote the cardinality
of $A$.
Consider the measure on $\Cpx^n$,
\[
\mut=\frac1{|A|}\sum_{\pi\in A}\mu_\pi.
\]
Let $\dt(1),\dt(2),\ldots,\dt(n)\in\LEu$ be random variables whose joint
distribution is $\mut$ and such that
\[
\{\dt(i)\mid1\le i\le n\}\qquad\mbox{and}\qquad\{t(i,j)\mid1\le i<j\le n\}
\]
are independent sets of random variables.
Let
\[
\Dwt=\sum_{i=1}^n\dt(i)\otimes e(i,i;n)\in\MEu_n
\]
and let $\Zwt=\Dwt+T(n)$.
Then $Z$ and $\Zwt$ have the same $*$--moments with respect to $\tau_n$.
\end{lemma}

\begin{proof}
We may introduce a uniformly distributed $A$--valued random variable, $\sigma$, that is
independent from $D$.
Then $\Dwt$ has the same distribution as the random matrix, $D_\sigma$, which at
a point $\omega\in\Omega$ takes the value
\[
D_\sigma(\omega)=\sum_{i=1}^nd(\sigma_{(\omega)}(i))_{(\omega)}\otimes e(i,i;n).
\]
Now each $*$--moment of $D_\sigma$ is the average over $\pi\in A$ of the
corresponding $*$--moments of the matrices $D_\pi$ described in~(\ref{eq:Dpi}).
By Lemma~\ref{lem:permind}, each of these is in turn equal to the
corresponding $*$--moment of $D$.
\end{proof}

Now we combine Theorem~\ref{thm:utgrm} and Lemma~\ref{lem-desymmetrize} to
obtain this section's main result.
\begin{thm}
\label{thm:*momPerm}
For each $n\in\Nats$ let $D_1(n),D_2(n)\in\MEu_n$ be diagonal random matrices, let
$T(n)\in\UTGRM(n,\frac1n)$ and let $Z_\iota(n)=D_\iota(n)+T(n)$,
($\iota=1,2$).
Suppose that $Z_1(n)$ converges in $*$--moments as $n\to\infty$.
For $\iota\in\{1,2\}$ and $n\in\Nats$ let $\nu_\iota^{(n)}$ be the measure on
$\Cpx^n$ that is the joint distribution of the $n$ diagonal entries of
$D_\iota(n)$, and let $\nut_\iota^{(n)}$ be the symmetrization of
$\nu_\iota^{(n)}$.
Suppose that for every $p\in\Nats$ and every
$r(1),s(1),\ldots,r(p),s(p)\in\Natsz$,
\begin{equation*}
\begin{split}
\lim_{n\to\infty}
\bigg(\int_{\Cpx^n}\lambda_1^{r(1)}\overline{\lambda_1}^{s(1)}
&\lambda_2^{r(2)}\overline{\lambda_2}^{s(2)}\cdots
\lambda_p^{r(p)}\overline{\lambda_p}^{s(p)}
\dif\nut_1^{(n)}(\lambda_1,\ldots,\lambda_n) \\
&-\int_{\Cpx^n}\lambda_1^{r(1)}\overline{\lambda_1}^{s(1)}
\lambda_2^{r(2)}\overline{\lambda_2}^{s(2)}\cdots
\lambda_p^{r(p)}\overline{\lambda_p}^{s(p)}
\dif\nut_2^{(n)}(\lambda_1,\ldots,\lambda_n)
\bigg)=0.
\end{split}
\end{equation*}
Then also $Z_2(n)$ converges in $*$--moments as $n\to\infty$, and its limit $*$--moments
are the same as those of $Z_1(n)$.
\end{thm}

\section{Asymptotically free random matrices}
\label{sec:asympfreerandmat}

This section concerns asymptotic freeness of self--adjoint
i.i.d.\ Gaussian random matrices $Y(t,n)$ and certain diagonal random matrices that are independent from the $Y(t,n)$,
this being a generalization of Voiculescu's pioneering
result~\cite{Voiculescu:RandMat}, which concerned constant diagonal matrices.
(See~\cite{Dy:FreeProdR} and~\cite{Voiculescu:StrRandMat} for some other generalizations.)
Just as, using a technique based on the polar decomposition, Voiculescu parlayed his asymptotic freeness results
for Gaussian random matrices into asymptotic freeness results for Haar distributed random unitary matrices,
so in this section do we prove asymptotic freeness of Haar distributed random unitary matrices $U(r,n)$ and certain random diagonal matrices
that are independent from the $U(r,n)$.
Finally, this section culminates in a result (Theorem~\ref{thm:Rdiag}) about matrix models for ($*$--free families of) $R$--diagonal elements.

See Notation~\ref{notation:UTGRM} for details of some terms used below.

\begin{notation}
\label{notation:GRM}
\renewcommand{\labelenumi}{(\roman{enumi})}
\begin{enumerate}
\item For a random matrix $X\in\MEu_n$, we will write $X\in\GRM(n,\sigma^2)$, (the acronym is for ``Gaussian random matrix'')
if the entries $x_{ij}$ of $X$ ($1\le i,j\le n$) satisfy that $(x_{ij})_{1\le i,j\le n}$ is an independent family
of random variables, each of which is complex $(0,\sigma^2)$--Gaussian.
Thus $X\in\GRM(n,\sigma^2)$ if and only if $X$ has density $(\pi\sigma^2)^{-n^2}\exp(-\frac1{\sigma^2}\Tr(X^*X))$ with respect to Lebesgue measure on $M_n(\Cpx)$.

\item Given $Y\in\MEu_n$ and $\sigma^2>0$, we will write $Y\in\SGRM(n,\sigma^2)$ (``self--adjoint Gaussian random matrix'') if the entries
$y_{ij}$ of $Y$ ($1\le i,j\le n$) satisfy
that $y_{ij}=\overline{y_{ji}}$ for all $i$ and $j$, that $y_{ij}$ is complex $(0,\sigma^2)$--Gaussian if $i\ne j$ and is real $(0,\sigma^2)$--Gaussian if $i=j$
and that $(y_{ij})_{1\le i\le j\le n}$ is an independent family of random variables.
Thus $Y\in\SGRM(n,\sigma^2)$ if and only if $Y$ has density $(\pi\sigma^2)^{-n^2/2}\exp(-\frac1{\sigma^2}\Tr(Y^*Y))$ with respect to Lebesgue measure on $M_n^{s.a.}$.

\item Given $U\in\MEu_n$, we will write $U\in\HURM(n)$, (``Haar unitary random matrix'') if $U$ is a random unitary matrix distributed according
to Haar measure on the $n\times n$ unitary matrices.
\end{enumerate}
\end{notation}

We begin with a preliminary result that is essentially just a combination of
Theorems~2.1 and~2.2 of~\cite{Voiculescu:RandMat}, in the case of random diagonal matrices.

\begin{lemma}
\label{lem:Voic2.1-2}
Let $S$ and $T$ be sets.
For any $n\in\Nats$ let $Y(s,n)\in\SGRM(n,1/n)$ ($s\in S$).
and consider diagonal random matrices $D(t,n)\in\MEu_n$ ($t\in T$).
Suppose that for some $t_0\in T$, $D(t_0,n)=I_n$ ($n\in\Nats$), that
$\{D(t,n)\mid t\in T\}$ is closed under multiplication ($n\in\Nats$)
and that $\{D(t,n)\mid t\in T\}$ converges in moments as $n\to\infty$.
Suppose that for every $n\in\Nats$ 
\begin{equation}
\label{eq:DYfam}
\bigg(\big\{D(t,n)\mid t\in T\big\},
\big(\{Y(s,n)\}\big)_{s\in S}\bigg)
\end{equation}
is an independent family of sets of matrix--valued random variables.
Finally, suppose that for every $m\in\Nats$, every $s_1,\ldots,s_m\in S$ and every $t_1,\ldots,t_m\in T$, the quantity
\begin{equation}
\label{eq:bddness}
\left|\tau_n\big(Y(s_1,n)D(t_1,n)\cdots Y(s_m,n)D(t_m,n)\big)\right|
\end{equation}
remains bounded as $n\to\infty$.
Then the following are equivalent:
\begin{enumerate}
\item
The family~\eqref{eq:DYfam} of sets of noncommutative random variables is asymptotically free as $n\to\infty$.
\item
\label{en-both}
Whenever $m\in\Nats$ is even, $t_1,\ldots,t_m\in T$ are fixed
and $\alpha:\{1,\ldots,m\}\to S$ is such that for every $s\in S$,
$\alpha^{-1}(s)$ has either two or zero elements,
\begin{enumerate}
\item
\label{en-alpha12}
if $\alpha(1)=\alpha(2)$ then
\begin{multline*}
\lim_{n\to\infty}\bigg(\tau_n\big(Y(\alpha(1),n)D(t_1,n)Y(\alpha(2),n)D(t_2,n)\cdots
Y(\alpha(m),n)D(t_m,n)\big) \\
-\tau_n\big(D(t_1,n)\big)\tau_n\big(D(t_2,n)Y(\alpha(3),n)D(t_3,n)\cdots
Y(\alpha(m),n)D(t_m,n)\big)\bigg)=0
\end{multline*}
\item
\label{en-alphaalt}
if $\alpha(p)\ne\alpha(p+1)$ for every $1\le p\le m-1$ and if
$\alpha(m)\ne\alpha(1)$ then
\[
\lim_{n\to\infty}\tau_n\big(Y(\alpha(1),n)D(t_1,n)Y(\alpha(2),n)D(t_2,n)\cdots
Y(\alpha(m),n)D(t_m,n)\big)=0.
\]
\end{enumerate}
\end{enumerate}
\end{lemma}
\begin{proof}
We may without loss of generality suppose $S=\Nats$ and $T=\Nats$.
That (1)$\implies$(2) follows from the last paragraph
of~\cite[2.1]{Voiculescu:RandMat} and the fact that the limit
moments of each $Y(s,n)$ are those of a centered semicircle law with second moment $1$.

To show (2)$\implies$(1) we will use~\cite[2.1]{Voiculescu:RandMat} and an
idea from the proof of~\cite[2.2]{Voiculescu:RandMat}.
Let $\omega$ be a nontrivial ultrafilter on $\Nats$.
On the algebra, $\Cpx\langle (T_s)_{s\in\Nats},(A_t)_{t\in\Nats}\rangle$, of
polynomials in noncommuting variables $(T_s)_{s\in\Nats}$ and
$(A_t)_{t\in\Nats}$, let $\phi_\omega$ be the tracial linear functional defined
by $\phi_\omega(p)=\lim_{n\to\omega}\tau_n(\pi_n(p))$,
where $\pi_n:\Cpx\langle (T_s)_{s\in\Nats},(A_t)_{t\in\Nats}\rangle\to\MEu_n$
is the algebra homomorphism defined by $\pi_n(T_s)=Y(s,n)$ and
$\pi_n(A_t)=D(t,n)$.
Let $\Delta$ denote the subalgebra of
$\Cpx\langle (T_s)_{s\in\Nats},(A_t)_{t\in\Nats}\rangle$
generated by $1$ and $\{A_t\mid t\in\Nats\}$.
We will check that the conditions 1$^\circ$ and 2$^\circ$
of~\cite[2.1]{Voiculescu:RandMat} hold for the sequence $(T_s)_{s\in\Nats}$
and the subalgebra $\Delta$ with respect to $\phi_\omega$.
Note that every element of $\Delta$ is a linear combination of words of the
form $A_{t_1}A_{t_2}\cdots A_{t_k}$ and that $\pi_n(A_{t_1}A_{t_2}\cdots
A_{t_k})=D(t,n)$ for some $t\in\Nats$.
Moreover, $\pi_n(1)=D(1,n)$.
Therefore $\pi_n(\Delta)=\lspan\{D(t,n)\mid t\in\Nats\}$ and hence
condition~1$^\circ$ of~\cite[2.1]{Voiculescu:RandMat} follows from the
boundedness of~(\ref{eq:bddness}) as $n\to\infty$.

To see that condition~2$^\circ$a of~\cite[2.1]{Voiculescu:RandMat} holds, it
suffices to see that if $m\in\Nats$ and if $\alpha:\{1,\ldots,m\}\to\Nats$ is
such that $\alpha^{-1}(\{\alpha(1)\})$ has only one element and if
$t_1,\ldots,t_m\in\Nats$ then
\begin{equation}
\label{eq:taunis0}
\forall n\in\Nats\qquad
\tau_n\big(Y(\alpha(1),n)D(t_1,n)Y(\alpha(2),n)D(t_2,n)\cdots
Y(\alpha(m),n)D(t_m,n)\big)=0.
\end{equation}
But this follows from the independence of $Y(\alpha(1),n)$ from all the other matrices
appearing in~(\ref{eq:taunis0}) and the fact that all
entries of $Y(\alpha(1),n)$ have expectation zero.

Now conditions~2$^\circ$b and~2$^\circ$c of~\cite[2.1]{Voiculescu:RandMat}
follow from the hypotheses~(\ref{en-alpha12}) and~(\ref{en-alphaalt}).
Therefore, by~\cite[2.1]{Voiculescu:RandMat}, given an injection
$\beta:\Nats\times\Nats\to\Nats$ and defining
\[
X_{m,k}=k^{-1/2}\sum_{j=1}^kT_{\beta(m,j)},
\]
the family of sets of noncommutative random variables,
\[
\bigg(\Delta,\big(\{X_{m,k}\}\big)_{m=1}^\infty\bigg)
\]
is asymptotically free with respect to $\phi_\omega$ as $k\to\infty$.
However, using the Gaussianity of the entries of the $Y(s,n)$, we see
that for every $k\in\Nats$, $\big(\Delta,(\{X_{m,k}\})_{m=1}^\infty\big)$ has
the same moments as $\big(\Delta,(\{T_m\})_{m=1}^\infty\big)$.
Hence $\big(\Delta,(\{T_m\})_{m=1}^\infty\big)$ is free with respect to
$\phi_\omega$.
Since $\omega$ was arbitrary, and since each $Y(s,n)$ converges in
moments as $n\to\infty$, this implies that
\[
\bigg(\{D(t,n)\mid t\in\Nats\},\big(\{Y(s,n)\}\big)_{s\in\Nats}\bigg)
\]
is asymptotically free as $n\to\infty$.
\end{proof}

\begin{thm}
\label{thm:indepdiag}
Let $S$ and $T$ be sets.
For $s\in S$ and $n\in\Nats$ let $Y(s,n)\in\SGRM(n,\frac1n)$.
For $t\in T$ and $n\in\Nats$ let
\[
D(t,n)=\sum_{i=1}^nd(i;n,t)\otimes e(i,i;n)\in\MEu_n
\]
be a diagonal random matrix, and suppose
that for some $t$ and every $n$, $D(1,n)=I_n$, that
$\{D(t,n)\mid t\in T\}$ is closed under multiplication
and that $\{D(t,n)\mid t\in T\}$ converges in moments as $n\to\infty$.
Suppose that for every $n\in\Nats$ 
\begin{equation*}
\bigg(\big\{D(t,n)\mid t\in T\big\},
\big(\{Y(s,n)\}\big)_{s\in S}\bigg)
\end{equation*}
is an independent family of sets of matrix--valued random variables.
Assume further that
\renewcommand{\labelenumi}{(\roman{enumi})}
\begin{enumerate}
\item
for every $t\in T$ and every $1\le p<\infty$,
\begin{equation*}
\sup_{\substack{n\in\Nats\\1\le i\le n}}\nm{d(i;n,t)}_p<\infty;
\end{equation*}
\item
for every $m,n\in\Nats$, $m\le n$, every $t_1,\ldots,t_m\in T$ and every permutation,
$\sigma$, of $\{1,\ldots,n\}$, the joint distribution of
\begin{equation*}
(d(1;n,t_1),d(2;n,t_2),\ldots,d(m;n,t_m))
\end{equation*}
is equal to the joint distribution of
\begin{equation*}
(d(\sigma(1);n,t_1),d(\sigma(2);n,t_2),\ldots,d(\sigma(m);n,t_m));
\end{equation*}
\item
for every $p\in\Nats$, every $t_1,\ldots, t_p\in T$ and every $p$--tuple,
$(i_1,\ldots,i_p)$, of distinct, positive integers, we have
\begin{equation*}
\lim_{n\to\infty}\bigg(E\big(d(i_1;n,t_1)d(i_2;n,t_2)\cdots d(i_p;n,t_p)\big)
-\prod_{j=1}^pE(d(i_j;n,t_j))\bigg)=0.
\end{equation*}
\end{enumerate}
Then the family
\[
\bigg(\{D(t,n)\mid t\in T\},\big(\{Y(s,n)\}\big)_{s\in S}\bigg)
\]
of sets of random variables converges in moments and is asymptotically free as $n\to\infty$.
\end{thm}

\begin{proof}
We will apply Lemma~\ref{lem:Voic2.1-2}.
Let us first show that the quantity~(\ref{eq:bddness}) remains bounded as
$n\to\infty$.
Write $a(i,j;n,s)$ for the $(i,j)$th entry of $Y(s,n)$.
We have
\begin{equation}
\begin{aligned}
\label{eq:taunYD}
\tau_n\big(Y(s_1,n)D(t_1,n)\cdots Y(s_m,n)&D(t_m,n)\big)=\\
=n^{-1}\sum_{i_1,\ldots,i_m\in\{1,\ldots,n\}}
&E\big(d(i_1;n,t_1)d(i_2;n,t_2)\cdots d(i_m;n,t_m)\big)\cdot \\
&\cdot E\big(a(i_0,i_1;n,s_1)a(i_1,i_2;n,s_2)\cdots a(i_{m-1},i_m;n,s_m)\big),
\end{aligned}
\end{equation}
where $i_0=i_m$.
Using the generalized H\"older inequality, we have
\begin{equation}
\label{eq:dHolder}
|E(d(i_1;n,t_1)\cdots d(i_m;n,t_m))|
\le\nm{d(i_1;n,t_1)\cdots d(i_m;n,t_m)}_1
\le\prod_{j=1}^m\nm{d(i_j;n,t_j)}_m.
\end{equation}
But by the assumption~(i), there is $c_2>0$ such that
\begin{equation}
\label{eq:bddmnroms}
\forall n\in\Nats\quad\forall i_1,\ldots,i_m\in\{1,\ldots,n\},\quad
\prod_{j=1}^m\nm{d(i_j;n,t_j)}_m\le c_2.
\end{equation}
Now consider
\begin{equation}
\label{eq:Eas}
E\big(a(i_0,i_1;n,s_1)a(i_1,i_2;n,s_2)\cdots a(i_{m-1},i_m;n,s_m)\big).
\end{equation}
From the nature of the entries $a(i,j;n,s)$, we see that the
quantity~(\ref{eq:Eas}) can be nonzero only if there is a bijection
$\gamma:\{1,\ldots,m\}\to\{1,\ldots,m\}$ such that $\gamma^2=\id$, $\gamma$ has
no fixed points and
\[
\forall j\in\{1,\ldots,m\},\qquad s_j=s_{\gamma(j)},\quad i_j=i_{\gamma(j)-1},
\quad i_{j-1}=i_{\gamma(j)}.
\]
One also calculates
\begin{align*}
\big|E\big(a(i_0,i_1;n,s_1)a(i_1,i_2;n,s_2)\cdots
a(i_{m-1},i_m;n,s_m)\big)\big|
&\le\prod_{j=1}^m\nm{a(i_{j-1},i_j;n,s_j)}_m \\
&\le n^{-m/2}\left(\tfrac m2\right)!
\end{align*}
Let us call a bijection, $\gamma$, of $\{1,\ldots,m\}$ without fixed points and
such that $\gamma^2=\id$, a {\it pairing} of $\{1,\ldots, m\}$ and for every
pairing $\gamma$ let
\[
\Theta(\gamma)=\big\{(i_1,\ldots,i_m)\in\{1,\ldots,n\}^m\mid
\forall j\in\{1,\ldots,m\},\,i_j=i_{\gamma(j)-1},\,
i_{j-1}=i_{\gamma(j)}\big\}.
\]
From the above estimates we obtain
\begin{equation}
\label{eq:taunYDsumgamma}
\big|\tau_n\big(Y(s_1,n)D(t_1,n)\cdots Y(s_m,n)D(t_m,n)\big)\big|
\le c_2\left(\tfrac m2\right)!\,n^{-(\frac m2+1)}\sum_\gamma|\Theta(\gamma)|,
\end{equation}
where the sum is over all pairings, $\gamma$ of $\{1,\ldots,m\}$.
To each pairing we associate the quotient graph, $G_\gamma$, of the clockwise oriented $m$--gon
graph obtained
by identifying with opposite orientation each pair of $j$th and $\gamma(j)$th edges.
The resulting graph has $m/2$ edges, hence at most $\frac m2+1$ vertices.
Consequently $|\Theta(\gamma)|\le n^{\frac m2+1}$, and the quantity
in~(\ref{eq:taunYDsumgamma}) is bounded by $c_2\left(\tfrac m2\right)!$ times the number of pairings,
which is finite and independent of $n$.
This shows that the quantities~(\ref{eq:bddness}) remain bounded as
$n\to\infty$.

We now show that~\ref{en-alpha12} and~\ref{en-alphaalt} of Lemma~\ref{lem:Voic2.1-2} are satisfied.
Let $\alpha$ be as described there.
Then
\begin{equation}
\begin{aligned}
\label{eq:taunYDalpha}
\tau_n\big(Y(\alpha(1),n)D(t_1,n)\cdots Y(\alpha&(m),n)D(t_m,n)\big)=\\
=n^{-1}\sum_{i_1,\ldots,i_m\in\{1,\ldots,n\}}E\big(&d(i_1;n,t_1)d(i_2;n,t_2)\cdots d(i_m;n,t_m)\big)\cdot \\
\cdot E\big(&a(i_0,i_1;n,\alpha(1))a(i_1,i_2;n,\alpha(2))\cdots a(i_{m-1},i_m;n,\alpha(m))\big),
\end{aligned}
\end{equation}
where we let $i_0=i_m$.
Consider the clockwise oriented $m$--gon graph, label the edges consecutively $e_1,e_2,\ldots,e_m$ and the vertices
$v_1,v_2,\ldots,v_m$ so that the vertices of the edge $e_j$ are $v_{j-1}$ and $v_j$, (mod $m$).
Let $G$ be the quotient of the $m$--gon graph obtained by identifying edges $j$ and $\alpha(j)$ with opposite orientation,
($1\le j\le m$).
The resulting identification of vertices of the $m$--gon graph gives an equivalence relation $\sim$
on $\{v_1,\ldots,v_m\}$ whose equivalence classes $F_1,F_2,\ldots,F_{k(G)}$ are precisely the lists of indices labeling the $k(G)$ vertices of $G$.
The expression
\begin{equation}
\label{eq:Ea}
E\big(a(i_0,i_1;n,\alpha(1))a(i_1,i_2;n,\alpha(2))\cdots a(i_{m-1},i_m;n,\alpha(m))\big)
\end{equation}
in~\eqref{eq:taunYDalpha} is nonzero if and only if whenever $1\le p,q\le m$ and $v_p\sim v_q$ then $i_p=i_q$, and then the value of~\eqref{eq:Ea}
is $n^{-m/2}$.
For each equivalence class $F_j=\{v_{p(1)},v_{p(2)},\ldots,v_{p(r_j)}\}$ of $\sim$ there is $t'_j\in T$ so that
$D(t'_j,n)=D(t_{p(1)},n)D(t_{p(2)},n)\cdots D(t_{p(r_j)},n)$;
for $p\in\{1,\ldots,n\}$ let $d(p;n,t'_j)$ be the $p$th diagonal entry of $D(t'_j,n)$.
Thus
\begin{equation}
\label{eq:taunYDalphaFj}
\begin{aligned}
\tau_n\big(Y(\alpha&(1),n)D(t_1,n)\cdots Y(\alpha(m),n)D(t_m,n)\big)=\\
&=n^{-(\frac m2+1)}\sum_{p_1,\ldots,p_{k(G)}\in\{1,\ldots,n\}}
E\big(d(p_1;n,t'_1)d(p_2;n,t'_2)\cdots d(p_{k(G)};n,t'_{k(G)})\big).
\end{aligned}
\end{equation}
Using the H\"older inequality estimate~\eqref{eq:dHolder} and~\eqref{eq:bddmnroms}, we see that the terms
\begin{equation*}
E\big(d(p_1;n,t'_1)d(p_2;n,t'_2)\cdots d(p_{k(G)};n,t'_{k(G)})\big)
\end{equation*}
in~\eqref{eq:taunYDalphaFj} are uniformly bounded in modulus.
Moreover, because $G$ has $m/2$ edges, it follows that $k(G)\le\frac m2+1$.
If $\alpha$ satisfies the hypothesis in~\ref{en-alphaalt} of Lemma~\ref{lem:Voic2.1-2}, then every vertex of the $m$--gon
graph is equivalent to at least one other vertex, so $k(G)\le m/2$ and the quantity~\eqref{eq:taunYDalphaFj} tends
to zero as $n\to\infty$, as required.
We have proved that~\ref{en-alphaalt} of Lemma~\ref{lem:Voic2.1-2} is satisfied.
(In fact, a similar analysis shows that the limit moment is zero unless the pairing of $\{1,\ldots,m\}$ given by $\alpha$
is non--crossing --- this was examined in a slightly different context in~\cite{Dy:FreeProdR}, but unfortunately
without the benefit of the idea of a non--crossing pairing.)

Suppose that $\alpha$ satisfies the hypothesis in~\ref{en-alpha12} of Lemma~\ref{lem:Voic2.1-2}, namely that $\alpha(1)=\alpha(2)$.
We are interested in the limit of the moment~\eqref{eq:taunYDalphaFj} as $n\to\infty$.
As the number of terms in the sum~\eqref{eq:taunYDalphaFj} where $p_i=p_j$ for some $i\ne j$ becomes
negligably small compared to $n^{\frac m2+1}$ as $n\to\infty$, we may in~\eqref{eq:taunYDalphaFj}
sum over only all distinct choices of $p_1,\ldots,p_{k(G)}\in\{1,\ldots,n\}$.
The assumption $\alpha(1)=\alpha(2)$ implies that $v_1$ is not equivalent to any other vertex under~$\sim$.
Therefore, renumbering if necessary, we may take $F_1=\{v_1\}$ and hence $t_1'=t_1$.
By hypotheses~(ii) and~(iii), we have that
\begin{align*}
\delta_n\eqdef E\big(d(p_1;n,t'_1)&d(p_2;n,t'_2)\cdots d(p_{k(G)};n,t'_{k(G)})\big) \\
&-E\big(d(p_1;n,t_1)\big)E\big(d(p_2;n,t'_2)\cdots d(p_{k(G)};n,t'_{k(G)})\big)
\end{align*}
is independent of the choice of distinct $p_1,\ldots,p_{k(G)}\in\{1,\ldots,n\}$ and that $\delta_n\to0$ as $n\to\infty$.
Moreover, an analysis of the quotient graph $G$ similar to that used to obtain~\eqref{eq:taunYDalphaFj}, and
keeping the same notation as in~\eqref{eq:taunYDalphaFj}, shows that
\begin{equation*}
\begin{aligned}
\tau_n\big(D(t_2,n)Y(&\alpha(3),n)D(t_3,n)\cdots Y(\alpha(m),n)D(t_m,n)\big)= \\
&=n^{-\frac m2}\sum_{p_2,\ldots,p_{k(G)}\in\{1,\ldots,n\}}
E\big(d(p_2;n,t'_2)\cdots d(p_{k(G)};n,t'_{k(G)})\big).
\end{aligned}
\end{equation*}
But then
\begin{align*}
&n^{-(\frac m2+1)}\sum_{p_1,\ldots,p_{k(G)}\in\{1,\ldots,n\}}
E\big(d(p_1;n,t_1)\big)E\big(d(p_2;n,t'_2)\cdots d(p_{k(G)};n,t'_{k(G)})\big) \\
&\;=\left(n^{-1}\sum_{p_1=1}^nE\big(d(p_1;n,t_1)\big)\right)
\left(n^{-\frac m2}\sum_{p_2,\ldots,p_{k(G)}\in\{1,\ldots,n\}}
E\big(d(p_2;n,t'_2)\cdots d(p_{k(G)};n,t'_{k(G)})\big)\right) \\
&\;=\tau_n\big(D(t_1,n)\big)\tau_n\big(D(t_2,n)Y(\alpha(3),n)D(t_3,n)\cdots Y(\alpha(m),n)D(t_m,n)\big).
\end{align*}
Taking the limit as $n\to\infty$, we can at will require $p_j$'s to be distinct and then relax this requirement;
using that $\delta_n\to0$, we obtain the conclusion of~\ref{en-alpha12} of Lemma~\ref{lem:Voic2.1-2}.
\end{proof}

Following Voiculescu's proof~\cite[Theorem 3.8]{Voiculescu:RandMat}, we will use polar decomposition to extend
the asymptotic freeness result of Theorem~\ref{thm:indepdiag} to include also Haar distributed random
unitary matrices.

\begin{thm}
\label{thm:HUdiag}
Let $R$, $S$ and $T$ be sets.
For every $n\in\Nats$ and $s\in S$ let $Z(s,n)\in\GRM(n,\frac1n)$ and for every $r\in R$ let $U(r,n)\in\HURM(n)$;
furthermore, for every $t\in T$ let
\begin{equation*}
D(t,n)=\sum_{i=1}^nd(i;n,t)\otimes e(i,i;n)\in\MEu_n
\end{equation*}
be diagnoal random matrices such that the family $(D(t,n))_{t\in T}$ is closed under multiplication and converges in moments as $n\to\infty$;
assume further that the entries $d(i;n,t)$ satisfy the conditions~(i), (ii) and~(iii) of the statement of Theorem~\ref{thm:indepdiag}.
Suppose that
\begin{equation}
\label{eq:ZUDindep}
\Bigl(\bigl(\{Z(s,n)\}\bigr)_{s\in S},\bigl(\{U(r,n)\}\bigr)_{r\in R},\{D(t,n)\mid t\in T\}\Bigr)
\end{equation}
is an independent family of sets of matrix--valued random variables.
Then the family
\begin{equation}
\label{eq:ZUDasfree}
\Bigl(\bigl(\{Z(s,n)^*,Z(s,n)\}\bigr)_{s\in S},\bigl(\{U(r,n)^*,U(r,n)\}\bigr)_{r\in R},\{D(t,n)\mid t\in T\}\Bigr)
\end{equation}
of sets of noncommutative random variables converges in moments and is asymptotically free as $n\to\infty$.
Furthermore, each $Z(s,n)$ converges in $*$--moments to a circular element and each $U(r,n)$ converges in $*$--moments to a Haar unitary, as $n\to\infty$.
\end{thm}
\begin{proof}
We shall use Theorem~\ref{thm:indepdiag} and adapt the proof of~\cite[Theorem 3.8]{Voiculescu:RandMat} to our situation.
\begin{claim}
\label{claim:ZDfree}
The family
\begin{equation*}
\Bigl(\bigl(\{Z(s,n)^*,Z(s,n)\}\bigr)_{s\in S},\{D(t,n)\mid t\in T\}\Bigr)
\end{equation*}
of sets of random variables is asymptotically free as $n\to\infty$.
\end{claim}
\begin{proof}
With
\begin{align*}
\RealPart Z(s,n)&=(Z(s,n)+Z(s,n)^*)/2 \\
\ImagPart Z(s,n)&=(Z(s,n)-Z(s,n)^*)/2i,
\end{align*}
each of $\RealPart Z(s,n)$ and $\ImagPart Z(s,n)$ is in $\SGRM(n,\frac1{2n})$ and
\begin{equation}
\label{eq:ReImZD}
\Bigl(\bigl(\{\RealPart Z(s,n)\}\bigr)_{s\in S},\bigl(\{\ImagPart Z(s,n)\}\bigr)_{s\in S},\{D(t,n)\mid t\in T\}\Bigr)
\end{equation}
is an independent family of sets of matrix--valued random variables.
Thus, by Theorem~\ref{thm:indepdiag}, each $\RealPart Z(s,n)$ and each $\ImagPart Z(s,n)$ converges in moments
to a semicircular element and the family~\eqref{eq:ReImZD} is asymptotically free as $n\to\infty$.
This proves Claim~\ref{claim:ZDfree}.
\end{proof}

Now take $W(r,n)\in\GRM(n,1/n)$ so that
\begin{equation*}
\Bigl(\bigl(\{W(r,n)\}\bigr)_{r\in R},\bigl(\{Z(s,n)\}\bigr)_{s\in S},\{D(t,n)\mid t\in T\}\Bigr)
\end{equation*}
is an independent family of sets of matrix--valued random variables.
By Claim~\ref{claim:ZDfree},
\begin{equation}
\label{eq:WZDfree}
\Bigr(\bigr(\{W(r,n)^*,W(r,n)\}\bigr)_{r\in R},\bigr(\{Z(s,n)^*,Z(s,n)\}\bigr)_{s\in S},\{D(t,n)\mid t\in T\}\Bigr)
\end{equation}
is asymptotically free as $n\to\infty$.
If $W(r,n)=V\bigl(W(r,n)^*W(r,n)\bigr)^{1/2}$ is the polar decomposition of $W(r,n)$, then $V$ is almost everywhere a unitary,
and is distributed according to Haar measure on the group of $n\times n$ unitaries.
Therefore, letting $U(r,n)$ be the polar part, $V$, of $W(r,n)$, these random unitary matrices satisfy the hypotheses of the theorem.
We will follow the proof of~\cite[Theorem 3.8]{Voiculescu:RandMat} to show the asymptotic freeness of~\eqref{eq:ZUDasfree}.

For $\epsilon>0$ let $Y_\epsilon(r,n)=W(r,n)\bigl(\epsilon+W(r,n)^*W(r,n)\bigr)^{-1/2}$.
\begin{claim}
\label{claim:YZDfree}
For every $\epsilon>0$, the family
\begin{equation}
\label{eq:YZDfree}
\Bigr(\bigr(\{Y_\epsilon(r,n)^*,Y_\epsilon(r,n)\}\bigr)_{r\in R},\bigr(\{Z(s,n)^*,Z(s,n)\}\bigr)_{s\in S},\{D(t,n)\mid t\in T\}\Bigr)
\end{equation}
is asymptotically free as $n\to\infty$.
\end{claim}
\begin{proof}
Given $A\in\MEu_n$ and $1\le p<\infty$, let $|A|_p=\big(\tau_n(A^*A)^{p/2}\big)^{1/p}$;
moreover, let $|A|_\infty$ be the essential supremum of the operator norm of $A$ evaluated at points of the underlying probability space.
Let $q$ be a noncommutative monomial in $d=2a+2b+c$ variables (for nonegative integers $a,b,c$), with coefficient equal to $1$.
Let $0<\delta\le1$.
By Step~I of the proof of~\cite[3.8]{Voiculescu:RandMat}, letting $f$ be the function $f(t)=(\epsilon+t)^{-1/2}$, there is a polynomial $Q_\delta$
such that, letting
\begin{equation*}
A_\delta(r,n)=W(r,n)Q_\delta\bigl(W(r,n)^*W(r,n)\bigr),
\end{equation*}
we have
\begin{equation*}
\limsup_{n\to\infty}|A_\delta(r,n)-Y_\epsilon(r,n)|_d<\delta.
\end{equation*}
Because $|Y_\epsilon(r,n)|_d\le|Y_\epsilon(r,n)|_\infty\le1$, it follows that $|A_\delta(r,n)|_d<1+\delta$.

The assumption~(i) of Theorem~\ref{thm:indepdiag} on the entries of $D(t,n)$ implies that for all $p\ge1$ and for all $t\in T$,
$\sup_{n\in\Nats}|D(t,n)|_p<\infty$.
Moreover, the convergence in $*$--moments as $n\to\infty$ of $Z(s,n)$ implies that $\sup_{n\ge1}|Z(s,n)|_p<\infty$ whenvever $p$ is an even integer;
however, as $|Z(s,n)|_p$ is increasing in $p$, this holds for all $1\le p<\infty$.
Fix $r_1,\ldots,r_a\in R$, $s_1,\ldots,s_b \in S$, $t_1,\ldots,t_c\in T$ and let
\begin{align*}
R_1(n,\epsilon)&=q\Bigl(\bigl(Y_\epsilon(r_i,n)^*\bigr)_{i=1}^a,\bigl(Y_\epsilon(r_i,n)\bigr)_{i=1}^a,
\bigl(Z(s_i,n)^*\bigr)_{i=1}^b,\bigl(Z(s_i,n)\bigr)_{i=1}^b,\bigl(D(t_i,n)\bigr)_{i=1}^c\Bigr) \\
R_2(n,\epsilon,\delta)&=q\Bigl(\bigl(A_\delta(r_i,n)^*\bigr)_{i=1}^a,\bigl(A_\delta(r_i,n)\bigr)_{i=1}^a,
\bigl(Z(s_i,n)^*\bigr)_{i=1}^b,\bigl(Z(s_i,n)\bigr)_{i=1}^b,\bigl(D(t_i,n)\bigr)_{i=1}^c\Bigr).
\end{align*}
We may chose a constant $K$ indepent of $\delta$ and large enough so that
\begin{align}
\label{eq:Kbd1}
\forall i\in\{1,\ldots,b\}\qquad&\limsup_{n\to\infty}|Z(s_i,n)|_d<K \\
\label{eq:Kbd2}
\forall i\in\{1,\ldots,c\}\qquad&\limsup_{n\to\infty}|D(t_i,n)|_d<K
\end{align}
Using H\"older's inequality we find
\begin{equation*}
\limsup_{n\to\infty}|R_1(n,\epsilon)-R_2(n,\epsilon,\delta)|_1\le2aK^{2b+c}(1+\delta)^{2a-1}\delta,
\end{equation*}
and therefore
\begin{equation}
\label{eq:tR12}
\lim_{\delta\to0}\;\limsup_{n\to\infty}\bigl|\tau_n\bigl(R_1(n,\epsilon)\bigr)-\tau_n\bigl(R_2(n,\epsilon,\delta)\bigr)\bigr|=0.
\end{equation}
The asymptotic freeness of
\begin{equation*}
\Bigr(\bigr(\{A_\delta(r,n)^*,A_\delta(r,n)\}\bigr)_{r\in R},\bigr(\{Z(s,n)^*,Z(s,n)\}\bigr)_{s\in S},\{D(t,n)\mid t\in T\}\Bigr)
\end{equation*}
follows from that of~\eqref{eq:WZDfree};
this together with~\eqref{eq:tR12} implies the asymptotic freeness of~\eqref{eq:YZDfree},
and claim~\ref{claim:YZDfree} is proved.
\end{proof}

Step~III of the proof of~\cite[3.8]{Voiculescu:RandMat} shows that
for every $\theta>0$ there is $\epsilon_0>0$ such that
\begin{equation}
\label{eq:YUtheta}
\limsup_{n\to\infty}|Y_\epsilon(r,n)-U(r,n)|_d<\theta
\end{equation}
whenever $0<\epsilon\le\epsilon_0$.
Let again $q$ be a noncommutative monomial having coefficient equal to $1$ and with degree $d=2a+2b+c$, and let $r_1,\ldots,r_a\in R$, $s_1,\ldots,s_b\in S$, $t_1,\ldots,t_c\in T$.
Let
\begin{equation*}
R_3(n)=q\Bigl(\bigl(U(r_i,n)^*\bigr)_{i=1}^a,\bigl(U(r_i,n)\bigr)_{i=1}^a,
\bigl(Z(s_i,n)^*\bigr)_{i=1}^b,\bigl(Z(s_i,n)\bigr)_{i=1}^b,\bigl(D(t_i,n)\bigr)_{i=1}^c\Bigr)
\end{equation*}
Letting $K$ be a constant so that~\eqref{eq:Kbd1} and~\eqref{eq:Kbd2} hold, we easily see using~\eqref{eq:YUtheta} and H\"older's inequality that
if $0<\epsilon\le\epsilon_0$ then
\begin{equation*}
\limsup_{n\to\infty}|R_1(n,\epsilon)-R_3(n)|_1\le2aK^{2b+c}\theta.
\end{equation*}
Therefore
\begin{equation*}
\lim_{\epsilon\to0}\;\limsup_{n\to\infty}\bigl|\tau_n\bigl(R_1(n,\epsilon)\bigr)-\tau_n\bigl(R_3(n)\bigr)\bigr|=0.
\end{equation*}
This, together with Claim~\ref{claim:YZDfree} shows that the family~\eqref{eq:ZUDasfree} is asympototically free as $n\to\infty$ and finishes the proof of the theorem.
\end{proof}

The following sort of result is standard, but a proof is provided here for completeness.
\begin{lemma}
\label{lem:vBv}
Let $(A,\phi)$ be a W$^*$--noncommutative probability space, let $B$ be a unital subalgebra of $A$, let $S$ be a set and
for every $s\in S$ let $v_s\in A$ be a unitary with $\phi(v_s)=0$.
Suppose the family
\begin{equation}
\label{eq:Bv}
\Bigl(B,\bigl(\{v_s^*,v_s\}\bigr)_{s\in S}\Bigr)
\end{equation}
of $|S|+1$ sets of noncommutative random vairables is free.
Then the family
\begin{equation}
\label{eq:vBv}
(v_sBv_s^*)_{s\in S}
\end{equation}
of unital subalgebras of $A$ is free.
\end{lemma}
\begin{proof}
From $*$--freeness of $B$ and $v_s$ we get $\phi(v_sbv_s^*)=\phi(b)$ for every $b\in B$.
Let $n\in\Nats$ and let $s_1,\ldots,s_n\in S$ be so that $s_j\ne s_{j+1}$ for every $j\in\{1,\ldots,n-1\}$.
For every $j\in\{1,\ldots,n\}$ let $b_j\in B$ be such that $\phi(b_j)=0$.
In order for freeness of~\eqref{eq:vBv} to hold, it will suffice that
\begin{equation*}
\phi\bigl((v_{s_1}b_1v_{s_1}^*)(v_{s_2}b_2v_{s_2}^*)\cdots(v_{s_n}b_nv_{s_n}^*)\bigr)=0.
\end{equation*}
But the above equality follows directly from freeness of~\eqref{eq:Bv}.
\end{proof}

Now we apply the asymptotic freeness results proved previously in this section to give some matrix models for
($*$--free families of) $R$--diagonal elements.

\begin{thm}
\label{thm:Rdiag}
Let $S$ be a set and for every $s\in S$ let
$X_s(n)\in\MEu_n$ and let $\sigma_{s,n}$ be the symmetrized joint distribution of the eigenvalues of $\bigl(X_s^*(n)X_s(n)\bigr)^{1/2}$.
Given $p\in\{1,\ldots,n\}$, let $\sigma_{s,n}^{(p)}$ be the marginal distribution of $\sigma_{s,n}$ corresponding to $p$ of the variables.
Suppose that for a compactly supported measure $\rho_s$ on $\Reals_+$ and for every $p\in\Nats$, $\sigma_{s,n}^{(p)}$ converges in moments
as $n\to\infty$ to the product measure $\btimes_1^p\rho_s$.
Suppose also that for any non--random $n\times n$ unitary matrix $U$, the distributions of $UX_s(n)$ and of $X_s(n)$ are the same.
\renewcommand{\labelenumi}{(\roman{enumi})}
\begin{enumerate}
\item
Fix $s\in S$.
Then $X_s(n)$ converges in $*$--moments to an element $u_sh_s$ of a noncommutative probability space,
where $u_s$ is a Haar unitary, $h_s\ge0$, $h_s$ has the same moments
as the measure $\rho_s$ and where the pair $\bigl(\{u_s,u_s^*\},\{h_s\}\bigr)$ of sets of noncommutative random variables is free.
\item
Suppose in addition that
\begin{equation}
\label{eq:Xsn}
\bigl(X_s(n)\bigr)_{s\in S}
\end{equation}
is a mutually independent family of matrix--valued random variables
and that the joint $*$--moments of~\eqref{eq:Xsn}
are the same as the joint $*$--moments of
\begin{equation}
\label{eq:UXUsn}
\bigl(U_s^{(1)}X_s(n)U_s^{(2)}\bigr)_{s\in S}
\end{equation}
whenever $U_s^{(1)}$ and $U_s^{(2)}$ are non--random unitary matrices, ($s\in S$).
Then the family~\eqref{eq:Xsn} is asymptotically $*$--free as $n\to\infty$.
\end{enumerate}
\end{thm}
\begin{proof}
For brevity we shall prove parts~(i) and~(ii) simultaneously;
while proving~(i), we may from the outset assume that the stronger hypotheses of~(ii) hold, because if we require $S$ to be
a single element then they will in any case be satisfied.
We may write $X_s(n)=V_s(n)H_s(n)$ where $V_s(n)$ is a random unitary matrix and $H_s(n)=\bigl(X_s(n)^*X_s(n)\bigr)^{1/2}$.
For every $s\in S$ let $W_s(n)$ be a random unitary matrix so that $D_s(n)=W_s(n)^*H_s(n)W_s(n)$ is diagonal,
so that the joint distribution of the diagonal entries of $D_s(n)$ is invariant under all permutations of the $n$ variables
and so that
\begin{equation*}
\bigl(\{D_s(n),V_s(n),W_s(n)\}\bigr)_{s\in S}
\end{equation*}
is a mutually independent family of sets of matrix--valued random variables.
Let $U_s^{(1)}(n),U_s^{(2)}(n)\in\HURM(n)$ be so that
\begin{equation*}
\Bigl(\bigl(\{H_s(n),V_s(n),W_s(n)\}\bigr)_{s\in S},\bigl(\{U_s^{(1)}(n)\}\bigr)_{s\in S},\bigl(\{U_s^{(2)}(n)\}\bigr)_{s\in S}\Bigr)
\end{equation*}
is a mutually independent family of sets of matrix--valued random variables.
It follows from the hypotheses of~(ii) that
\begin{equation*}
\bigl(U_s^{(1)}(n)X_s(n)U_s^{(2)}(n)\bigr)_{s\in S}
\end{equation*}
 has the same joint $*$--moments as the family~\eqref{eq:Xsn}.
We have
\begin{equation*}
U_s^{(1)}(n)X_s(n)U_s^{(2)}(n)=\bigl(U_s^{(1)}(n)V_s(n)U_s^{(2)}(n)\bigr)\bigl(U_s^{(2)}(n)^*W_s(n)D_s(n)W_s(n)^*U_s^{(2)}(n)\bigr).
\end{equation*}
Let\begin{align*}
\Vt_s(n)&=U_s^{(1)}(n)V_s(n)U_s^{(2)}(n) \\
\Wt_s(n)&=U_s^{(2)}(n)^*W_s(n).
\end{align*}
Then $\Vt_s(n),\Wt_s(n)\in\HURM(n)$ and
\begin{equation}
\label{eq:DVW}
\Bigl(\bigl(D_s(n)\bigr)_{s\in S},\bigl(\Vt_s(n)\bigr)_{s\in S},\bigl(\Wt_s(n)\bigr)_{s\in S}\Bigr)
\end{equation}
is an independent family of matrix--valued random variables.

Let
\begin{equation*}
\Delta(n)=\{I_n\}\cup\{D_{s_1}(n)D_{s_2}(n)\cdots D_{s_q}(n)\mid q\in\Nats,\,s_1,\ldots,s_q\in S\}.
\end{equation*}
By hypothesis, each $D_s(n)$ converges in moments as $n\to\infty$.
Since $\Delta(n)$ forms a commuting family of self--adjoint random matrices,
and since the family
\begin{equation}
\label{eq:Ds}
\bigl(D_s(n)\bigr)_{s\in S,}
\end{equation}
is independent, it follows that $\Delta(n)$ converges in moments as $n\to\infty$;
moreover, the subfamily~\eqref{eq:Ds} converges in $*$--moments to a family $(d_s)_{s\in S}$ is some W$^*$--noncommutative probability space $(A,\phi)$,
where $d_s$ is positive and  has the same moments as the measure $\sigma_s$, and where for distinct $s_1,\ldots,s_m\in S$ and any $k_1,\ldots,k_m\in\Nats$,
\begin{equation}
\label{eq:phids}
\phi(d_{s_1}^{k_1}d_{s_2}^{k_2}\cdots d_{s_m}^{k_m})=\prod_{j=1}^m\phi(d_{s_j}^{k_j}).
\end{equation}
We shall show that the entries of the set $\Delta(n)$ of diagonal random matrices satisfy the properties~(i), (ii) and~(iii) in the statement
of Theorem~\ref{thm:indepdiag}.
Let $d_s(i,n)$ denote the $i$th diagonal entry of $D_s(n)$.
Note that $E(d_s(i,n)^k)$ stays bounded (in fact converges) as $n\to\infty$, for every $k\in\Nats$ and $s\in S$.
This, together with the independence of the family~\eqref{eq:Ds},
implies condition~(i).
Condition~(ii) follows from the independence of~\eqref{eq:Ds} and the fact that the joint distribution of the diagonal entries of each $D_s(n)$
is invariant under permutations of the $n$ variables.
Because $\sigma_{s,n}^{(p)}$ converges to a product measure, we have for every $s\in S$, $m\in\Nats$, $k_1,\ldots,k_m\in\Nats$
and every $p$--tuple $(i_1,i_2,\ldots,i_m)$ of distinct, positive integers, that
\begin{equation*}
\lim_{n\to\infty}\bigg(E\big(d_s(i_1,n)^{k_1}d_s(i_2,n)^{k_2}\cdots d_s(i_m,n)^{k_m}\big)
-\prod_{j=1}^mE(d_s(i_j,n)^{k_j})\bigg)=0.
\end{equation*}
This together with the independence of~\eqref{eq:Ds} implies condition~(iii).
Hence we conclude from Theorem~\ref{thm:HUdiag} that
\begin{equation*}
\Bigl(\{D_s(n)\mid s\in S\},\bigl(\{\Vt_s(n)^*,\Vt_s(n)\}\bigr)_{s\in S},\bigl(\{\Wt_s(n)^*,\Wt_s(n)\}\bigr)_{s\in S}\Bigr)
\end{equation*}
is asymptotically free as $n\to\infty$.
Therefore the family~\eqref{eq:DVW} converges in $*$--moments to a family
\begin{equation*}
\bigl((d_s)_{s\in S},(v_s)_{s\in S},(w_s)_{s\in S}\bigr)
\end{equation*}
in some W$^*$--noncommutative probability space, where the joint $*$--moments of $(d_s)_{s\in S}$ are as described above,
where each $v_s$ and each $w_s$ is a Haar unitary and where
\begin{equation}
\label{eq:dvw}
\Bigl(\{d_s\mid s\in S\},\bigl(\{v_s^*,v_s\}\bigr)_{s\in S},\bigl(\{w_s^*,w_s\}\bigr)_{s\in S}\Bigr)
\end{equation}
is free.
Therefore, the family~\eqref{eq:Xsn} converges in $*$--moments as $n\to\infty$ to the family
\begin{equation*}
\bigl(v_s(w_sd_sw_s^*)\bigr)_{s\in S}.
\end{equation*}
It is clear that $w_sd_sw_s^*$ has the same moments as $d_s$, namely the same moments as the measure $\sigma_s$.
From the freeness of~\eqref{eq:dvw} and Lemma~\ref{lem:vBv}, it follows that the family
\begin{equation*}
\bigl((v_s)_{s\in S},(w_sd_sw_s^*)_{s\in S}\bigr)
\end{equation*}
is $*$--free, and the theorem is proved.
\end{proof}

\section{Upper triangular representations of circular free Poisson elements}
\label{sec:utcfP}

In this section, the random matrix results of~\S\ref{sec:utrandmat} and~\S\ref{sec:asympfreerandmat} are used, together
with results of Dyson and others, to give upper triangular matrix models of circular free Poisson elements, and finally
to give an upper triangular realization of a circular free Poisson element.
An outline of the contents of this section is as follows:
a first intermediate goal is a unitarily invariant matrix model for a circular free Poisson element (Theorem~\ref{thm:RPmm});
next, a result of Dyson is quoted (Theorem~\ref{thm:Dyson}) and used to convert the unitarily invariant matrix model to an upper
triangular matrix model for a circular free Poisson element (Corollary~\ref{cor:Dysondet});
then the diagonal elements of this upper triangular matrix model are decoupled and desymmetrized so as to yield,
in the limit as matrix size increases without bound, a triangular
realization of a circular free Poisson element (Theorem~\ref{thm:utfrPsn}).

The following is due to Bronk~\cite{Bronk}; see also~\cite[\S5]{HT}.
\begin{thm}
\label{thm:Bronk}
Let $c\ge1$ and let $Y$ be an $n\times n$ random matrix whose density with respect to Lebesgue measure on $M_n(\Cpx)$ is
\begin{equation*}
K_{c,n}^{(1)}|\det Y|^{2(c-1)n}\exp\big(-n\Tr(Y^*Y)\big),
\end{equation*}
where $K_{c,n}^{(1)}$ is a constant.
Then the symmetrized joint distribution of the eigenvalues of $Y^*Y$ has density
\begin{equation}
\label{eq:HTi}
K^{(2)}_{c,n}\left(\prod_{i=1}^n\lambda_i\right)^{2(c-1)n}\left(\prod_{1\le i<j\le n}(\lambda_i-\lambda_j)^2\right)
\exp\Bigl(-n\sum_{i=1}^n\lambda_i\Bigr)
\end{equation}
with respect to Lebesgue measure on $(\Reals_+)^n$, where $K^{(2)}_{c,n}$ is a constant.
\end{thm}

The next theorem is a corollary of a result of Hewitt and Savage~\cite{HS}.
\begin{thm}
\label{thm:ProdMeas}
Let $(\Omega,\Eb)$ be a standard Borel space.
Let $\sigma$ be a Borel probability measure on the product set $\Omega^\Nats=\prod_{n=1}^\infty\Omega$
endowed with the product topology.
Let $\sigma_1$ and $\sigma_2$ be the probability measures on $\Omega$ and $\Omega\times\Omega$, respectively, determined by
\begin{align*}
\sigma_1(A)&=\sigma(A\times\Omega\times\Omega\times\cdots),\quad A\in\Eb \\
\sigma_2(A_1\times A_2)&=\sigma(A_1\times A_2\times\Omega\times\Omega\times\cdots),\quad A_1,A_2\in\Eb.
\end{align*}
Suppose that
\renewcommand{\labelenumi}{(\roman{enumi})}
\begin{enumerate}
\item $\sigma$ is invariant under all finite permutations of coordinates in $\Omega^\Nats$,
(i.e.\ those permutations leaving all but finitely many coordinates fixed);
\item $\sigma_2=\sigma_1\times\sigma_1$.
\end{enumerate}
Then $\sigma$ is equal to the product measure $\btimes_{n=1}^\infty\sigma_1$.
\end{thm}
\begin{proof}
Since any noncountable standard Borel space is Borel isomorphic to the unit interval, and since $(\Nats,2^\Nats)$ is
Borel isomorphic to the one--point compactification of $\Nats$, it is no loss of generality to assume
that $\Omega$ is a separable compact Hausdorff space and $\Eb$ is the Borel $\sigma$--algebra associated to this topology.

For any compact set $K$, let $P(K)$ denote the set of Borel probability measures on $K$.
Consider the folowing subsets of $P(\Omega^\Nats)$:
\begin{align*}
\Pt&=\{\btimes_{n=1}^\infty\mu\mid\mu\in P(\Omega)\}, \\
\St&=\{\nu\in P(\Omega^\Nats)\mid\nu\text{ is invariant under all finite permutations of the coordinates of }\Omega^\Nats\}.
\end{align*}
Clearly $\Pt\subseteq\St$.
By~\cite[Theorem 7.2]{HS}, every $\nu\in\St$ has a representation 
\begin{equation*}
\nu=\int_{P(\Omega)}\left(\btimes_{n=1}^\infty\mu\right)\dif\rho(\mu),
\end{equation*}
for a unique $\rho\in P(P(\Omega))$.
In fact, (see~\cite[Theorem 3.1]{Sto}), $\Pt$ is the set of extreme points of the compact simplex $\St$.

Now let $\sigma$ be as in the formulation of the theorem.
Using hypothesis~(i) we have
\begin{equation*}
\sigma=\int_{P(\Omega)}\left(\btimes_{n=1}^\infty\mu\right)\dif\rho(\mu)
\end{equation*}
for a unique $\rho\in P(P(\Omega))$.
In particular
\begin{equation*}
\sigma_1=\int_{P(\Omega)}\mu\,\dif\rho(\mu) \qquad\text{and}\qquad
\sigma_2=\int_{P(\Omega)}(\mu\times\mu)\,\dif\rho(\mu).
\end{equation*}
By the assumption on $\Omega$, the space $C(\Omega)$ of complex valued continuous functions on $\Omega$
is a separable Banach space (in the uniform norm), so we may let $F$ be a countable dense subset of $C(\Omega)$.
Given  $f\in C(\Omega)$ and $\lambda\in P(\Omega)$ let us write
\begin{equation*}
\lambda(f)=\int_\Omega f\dif\lambda.
\end{equation*}
With this notation, we have for all $f\in F$,
\begin{align*}
&\int_{P(\Omega)}|\mu(f)-\sigma_1(f)|^2\dif\rho(\mu) \\
&=\int_{P(\Omega)}(\mu\times\mu)(f\otimes\overline f)\dif\rho(\mu)
-2\RealPart\left(\overline{\sigma_1(f)}\int_{P(\Omega)}\mu(f)\dif\rho(\mu)\right)+|\sigma_1(f)|^2 \\
&=\sigma_2(f\otimes\overline f)-2\RealPart\bigl(\overline{\sigma_1(f)}\sigma_1(f)\bigr)+|\sigma_1(f)|^2 \\
&=\sigma_2(f\otimes\overline f)-|\sigma_1(f)|^2.
\end{align*}
But hypothesis~(ii) shows that the above quantity is zero.
Hence $\mu(f)=\sigma_1(f)$ for all $f\in F$, for $\rho$--almost all $\mu\in P(\Omega)$.
Hence $\mu=\sigma_1$ for $\rho$--almost all $\mu\in P(\Omega)$, which implies $\rho=\delta_{\sigma_1}$,
the Dirac measure at the point $\sigma_1$.
Therefore $\sigma=\btimes_{n=1}^\infty\sigma_1$.
\end{proof}

\begin{lemma}
\label{lem:prodmarg}
Let $c\ge1$.
Given $n\in\Nats$ let $\mu_n$ be the measure having density~\eqref{eq:HTi} with respect to Lebesgue measure on $\Reals_+^n$.
Given $p\in\{1,2,\ldots,n\}$ let $\mu_n^{(p)}$ be the marginal distribution of $\mu_n$ corresponding
to the variables $\lambda_1,\ldots,\lambda_p$.
Fix $p\in\Nats$.
Then the distribution $\mu_n^{(p)}$ converges in the weak$^*$ topology
as $n\to\infty$ to the product measure $\btimes_1^p\tau$, where $\tau$ has density with respect to Lebesgue measure
\begin{equation}
\label{eq:munp}
\frac{\dif\tau}{\dif\lambda}=\frac{\sqrt{(\lambda-a)(b-\lambda)}}{2\pi\lambda}1_{[a,b]}(\lambda),
\end{equation}
with $a=(1-\sqrt c)^2$ and $b=(1+\sqrt c)^2$.
Moreover, if $f$ is a continuous function on $[0,\infty)^p$ with polynomial growth, in the sense that
$f(t_1,\ldots,t_p)\le K^{(3)}(1+t_1^{k_1}t_2^{k_2}\cdots t_p^{k_p})$
for some constant $K^{(3)}>0$ and positive integers $k_1,\ldots,k_p$, then
\begin{equation}
\label{eq:fpolygth}
\lim_{n\to\infty}\int_{\Reals_+^p}f\dif\mu_n^{(p)}=\int_{\Reals_+^p}f\dif\bigl(\btimes_1^p\tau\bigr).
\end{equation}
\end{lemma}
\begin{proof}
It will be more convenient to consider the measure $\sigma_n$ whose density with respect to Lebesgue measure on $\Reals_+^n$ is
\begin{equation*}
K^{(4)}_{c,n}\left(\prod_{i=1}^n\lambda_i\right)^{2(c-1)n}\left(\prod_{1\le i<j\le n}(\lambda_i-\lambda_j)^2\right)
\exp\left(-\sum_{i=1}^n\lambda_i\right),
\end{equation*}
for some constant $K^{(4)}_{c,n}$;
thus $\sigma_n$ is the push forward measure of $\mu_n$ under the transformation
$(\lambda_1,\lambda_2,\ldots)\mapsto(n\lambda_1,n\lambda_2,\ldots)$.
We will find the limit as $n\to\infty$ of the marginal distributions, $\sigma_n^{(p)}$, of $\sigma_n$
corresponding to the variables $\lambda_1,\ldots,\lambda_p$.
Let $\alpha=2(c-1)n$ and let $\phi_0^{(\alpha)},\phi_1^{(\alpha)},\phi_2^{(\alpha)},\cdots$ be the polynomials
obtained via Gram--Schmidt orthonormalizaton of $1,\lambda,\lambda^2,\ldots$ 
in $L^2([0,\infty),\lambda^\alpha e^{-\lambda}\dif\lambda)$.
Thus $\phi_k^{(\alpha)}(\lambda)=\sqrt{\frac{n!}{\Gamma(\alpha+n+1)}}L^\alpha_k(\lambda)$, where $L^\alpha_k$ are the (generalized)
Laguerre polynomials.
Using the Vandermonde determinant we have
\begin{align*}
\prod_{1\le i<j\le n}(\lambda_j-\lambda_i)&=\det\left(
\begin{matrix}
1&1&\cdots&1\\ \lambda_1&\lambda_2&\cdots&\lambda_n\\ \vdots&\vdots&&\vdots\\ \lambda_1^{n-1}&\lambda_2^{n-1}&\cdots&\lambda_n^{n-1}
\end{matrix}\right) \\[2ex]
&=K^{(5)}_{\alpha,n}\det\left(
\begin{matrix}
\phi_0^{(\alpha)}(\lambda_1)&\phi_0^{(\alpha)}(\lambda_2)&\cdots&\phi_0^{(\alpha)}(\lambda_n) \\[1ex]
\phi_1^{(\alpha)}(\lambda_1)&\phi_1^{(\alpha)}(\lambda_2)&\cdots&\phi_1^{(\alpha)}(\lambda_n) \\
\vdots&\vdots&&\vdots\\
\phi_{n-1}^{(\alpha)}(\lambda_1)&\phi_{n-1}^{(\alpha)}(\lambda_2)&\cdots&\phi_{n-1}^{(\alpha)}(\lambda_n)
\end{matrix}\right)
\end{align*}
for some constant $K^{(5)}_{\alpha,n}$.
Therefore, the density of $\sigma_n$ with respect to Lebesgue measure on $\Reals_+^n$ is
\pagebreak[0]
\begin{equation*}
D_n(\lambda_1,\ldots,\lambda_n)=K^{(6)}_{c,n}\left(\prod_{i=1}^n\lambda_i\right)^\alpha
\left(\sum_{\pi\in S_n}\sign(\pi)\prod_{i=1}^n\phi_{\pi(i)-1}^{(\alpha)}(\lambda_i)\right)^2
\exp\left(-\sum_{i=1}^n\lambda_i\right),
\end{equation*}
for a constant $K^{(6)}_{c,n}$.
Writing
\begin{equation*}
\prod_{i=1}^n\phi_{\pi(i)-1}^{(\alpha)}(\lambda_i)
=\phi_{\pi(1)-1}^{(\alpha)}(\lambda_1)\otimes\phi_{\pi(2)-1}^{(\alpha)}(\lambda_2)\otdt\phi_{\pi(n)-1}^{(\alpha)}(\lambda_n)
\end{equation*}
and noting that as $\pi$ ranges over the permutation group $S_n$ these form an orthonormal family
with respect to the measure
$\bigl(\prod_{i=1}^n\lambda_i\bigr)^\alpha\exp\bigl(-\sum_{i=1}^n\lambda_i\bigr)\dif\lambda_1\cdots\dif\lambda_n$
on $\Reals_+^n$, we find $K_{c,n}^{(6)}=(n!)^{-1}$.
Moreover, the density with respect to Lebesgue measure on $\Reals_+$ of the marginal distribution $\sigma_n^{(1)}$ is
\begin{equation*}
D_{n,1}(\lambda_1)\eqdef\int_{\Reals_+^{n-1}}D_n(\lambda_1,\ldots,\lambda_n)\dif\lambda_2\cdots\dif\lambda_n
=\frac1n\Bigl(\sum_{k=0}^{n-1}\phi_k^{(\alpha)}(\lambda_1)^2\Bigr)\lambda_1^\alpha e^{-\lambda_1}.
\end{equation*}
But then the treatment in~\S6 of~\cite{HT} shows that $\mu_n^{(1)}$ converges in the weak$^*$ topology and in moments as $n\to\infty$ to
$\tau$.

The density with respect to Lebesgue measure on $\Reals_+$ of the marginal distribution $\sigma_n^{(2)}$ is
\begin{align*}
D_{n,2}(\lambda_1)&\eqdef\int_{\Reals_+^{n-2}}D_n(\lambda_1,\ldots,\lambda_n)\dif\lambda_3\cdots\dif\lambda_n \\
&=\begin{aligned}[t]
 \frac1{n(n-1)}&(\lambda_1\lambda_2)^\alpha e^{-(\lambda_1+\lambda_2)}\cdot \\
 \cdot&\sum_{\substack{0\le k,\ell\le n-1\\ k\ne\ell}}\phi_k^{(\alpha)}(\lambda_1)\phi_\ell^{(\alpha)}(\lambda_2)
 \Bigl(\phi_k^{(\alpha)}(\lambda_1)\phi_\ell^{(\alpha)}(\lambda_2)-\phi_\ell^{(\alpha)}(\lambda_1)\phi_k^{(\alpha)}(\lambda_2)\Bigr)
 \end{aligned} \\
&=\frac n{n-1}D_{n,1}(\lambda_1)D_{n,1}(\lambda_2)-
\frac1{n(n-1)}(\lambda_1\lambda_2)^\alpha e^{-(\lambda_1+\lambda_2)}
\left(\sum_{j=0}^{n-1}\phi_j^{(\alpha)}(\lambda_1)\phi_j^{(\alpha)}(\lambda_2)\right)^2_.
\end{align*}
As elements of $C_0(\Reals_+^2)^*$, we thus have
\begin{equation}
\label{eq:sig21}
\begin{aligned}
\nm{\sigma_n^{(2)}-&\sigma_n^{(1)}\otimes\sigma_n^{(1)}}\le \\
&\le\frac1{n-1}\int_0^\infty\int_0^\infty D_{n,1}(\lambda_1)D_{n,1}(\lambda_2)\dif\lambda_1\dif\lambda_2\;+ \\
&\;\;+\frac1{n(n-1)}\sum_{j=0}^{n-1}\int_0^\infty\int_0^\infty\phi_j^{(\alpha)}(\lambda_1)^2\phi_j^{(\alpha)}(\lambda_2)^2
(\lambda_1\lambda_2)^\alpha e^{-(\lambda_1+\lambda_2)}\dif\lambda_1\dif\lambda_2 \\
&=\frac2{n-1}.
\end{aligned}
\end{equation}
Since we know that $\mu_n^{(1)}$ converges in the weak$^*$ topology as $n\to\infty$ to
$\tau$, it follows from~\eqref{eq:sig21} that
$\mu_n^{(2)}$ converges in weak$^*$ topology as $n\to\infty$ to $\tau\times\tau$.

Consider the measures $\mut_n=\mu_n\times\delta_0\times\delta_0\times\cdots$ on $[0,\infty)^\Nats$
and let $\nu$ be a w$^*$ cluster point in $C_0([0,\infty)^\Nats)^*$ of these.
Let $\nu^{(p)}$ be the marginal distribution of $\nu$ corresponding to the first $p$ coordinates of $[0,\infty)^\Nats$.
Then from what we have proved above we have
\renewcommand{\labelenumi}{$\bullet$}
\begin{enumerate}
\item $\nu$ is invariant under finite permutations of the coordinates in $[0,\infty)^\Nats$;
\item $\nu^{(1)}=\tau$;
\item $\nu^{(2)}=\tau\times\tau$.
\end{enumerate}
Hence, by Theorem~\ref{thm:ProdMeas}, $\nu=\displaystyle\btimes_{p=1}^\infty\nu^{(1)}$.
Since $\nu$ was an arbitrary cluster point of $(\mut_n)_{n=1}^\infty$ it follows that $\mut_n$ converges in weak$^*$ topology
to $\displaystyle\btimes_{p=1}^\infty\tau$ as $n\to\infty$.
Therefore, for all $p\in\Nats$, the marginal distribution $\nu_n^{(p)}$ converges in weak$^*$ topology to the measure $\btimes_1^p\tau$
as $n\to\infty$.

It remains to show that~\eqref{eq:fpolygth} holds whenever $f$ is of polynomial growth.
\begin{claim}
\label{claim:hpos}
Let $p\in\Nats$ and let $h$ be a positive continuous function on $[0,\infty)^p$.
Then
$\liminf_{n\to\infty}\int h\dif\mu_n^{(p)}\ge\int h\dif\nu^{(p)}$.
\end{claim}
\begin{proof}
Choose $h_j\in C_0\bigl([0,\infty)^p\bigr)$,
$h_j\ge0$, so that $h_j$ increases pointwise to $h$ as $j\to\infty$.
Then for all $j\ge1$,
\begin{equation*}
\liminf_{n\to\infty}\int h\dif\mu_n^{(p)}\ge\liminf_{n\to\infty}\int h_j\dif\mu_n^{(p)}
=\int h_j\dif\nu^{(p)}.
\end{equation*}
But $\nu^{(p)}=\btimes_1^p\nu^{(1)}$ is supported on $[a,b]^p$;
therefore $\sup_j\int h_j\dif\nu^{(p)}=\int h\dif\nu^{(p)}$, and the claim is proved.
\end{proof}

\begin{claim}
\label{claim:mgfg}
Let $p\in\Nats$ and suppose $f$ and $g$ are continuous functions on $[0,\infty)^p$ satisfying $g\ge0$ and $-g\le f\le g$,
and suppose that $\lim_{n\to\infty}\int g\dif\mu_n^{(p)}=\int g\dif\nu^{(p)}$.
Then $\lim_{n\to\infty}\int f\dif\mu_n^{(p)}=\int f\dif\nu^{(p)}$.
\end{claim}
\begin{proof}
Applying Claim~\ref{claim:hpos} to $g-f$ gives
$\limsup_{n\to\infty}\int f\dif\mu_n^{(p)}\le\int f\dif\nu^{(p)}$,
while applying Claim~\ref{claim:hpos} to $g+f$ yields
$\liminf_{n\to\infty}\int f\dif\mu_n^{(p)}\ge\int f\dif\nu^{(p)}$;
the claim is proved.
\end{proof}

In order to finish the proof of the lemma, it will suffice to show
\begin{equation}
\label{eq:munmom}
\lim_{n\to\infty}\int\lambda_1^{k_1}\lambda_2^{k_2}\cdots\lambda_p^{k_p}\dif\mu_n^{(p)}
=\int\lambda_1^{k_1}\lambda_2^{k_2}\cdots\lambda_p^{k_p}\dif\nu^{(p)},
\end{equation}
for every $p\in\Nats$ and all integers $k_1,k_2,\ldots,k_p\ge0$.
Letting $\ell=k_1+k_1+\cdots+k_p$, we have
$0\le\lambda_1^{k_1}\lambda_2^{k_2}\cdots\lambda_p^{k_p}\le1+\lambda_1^\ell+\lambda_2^\ell+\cdots+\lambda_p^\ell$.
Moreover, because $\mu_n^{(1)}$ converges in moments to $\nu^{(1)}$, we have
\begin{align*}
\lim_{n\to\infty}\int(1+\lambda_1^\ell+\lambda_2^\ell+\cdots+\lambda_p^\ell)\dif\mu_n^{(p)}
&=1+p\lim_{n\to\infty}\int\lambda_1^\ell\dif\mu_n^{(1)}
=1+p\int\lambda_1^\ell\dif\nu^{(1)}= \\
&=\int(1+\lambda_1^\ell+\lambda_2^\ell+\cdots+\lambda_p^\ell)\dif\nu^{(p)}.
\end{align*}
Now~\eqref{eq:munmom} follows from Claim~\ref{claim:mgfg}, and the lemma is proved.
\end{proof}

\begin{thm}
\label{thm:RPmm}
Let $c\ge1$ and let $Y(n)$ be an $n\times n$ random matrix whose density with respect to Lebesgue measure on $M_n(\Cpx)$ is
\begin{equation*}
K^{(1)}_{c,n}|\det Y|^{2(c-1)n}\exp\bigl(-n\Tr(Y^*Y)\bigr).
\end{equation*}
Then $Y(n)$ converges in $*$--moments as $n\to\infty$ to a circular free Poisson element of parameter $c$.
\end{thm}
\begin{proof}
Clearly for every non--random $n\times n$ unitary matrix $U$, the distribution of $UY(n)$ is equal to the distribution of $Y(n)$.
Let $\sigma_n$ be the symmetrized joint distribution of the eigenvalues of $\bigl(Y(n)^*Y(n)\bigr)^{1/2}$ and let $\mu_n$ be
the symmetrized joint distribution of the eigenvalues of $Y(n)^*Y(n)$.
For $p\in\{1,\ldots,n\}$ let $\sigma_n^{(p)}$, respectively $\mu_n^{(p)}$, be the marginal distribution of $\sigma_n$, respectively $\mu_n$,
corresponding to the first $p$ variables.
Given $k_1,\ldots,k_p\in\Nats\cup\{0\}$,
\begin{equation*}
\int\lambda_1^{k_1}\lambda_2^{k_2}\cdots\lambda_p^{k_p}\dif\sigma_n^{(p)}(\lambda_1,\ldots,\lambda_p)=
\int\lambda_1^{k_1/2}\lambda_2^{k_2/2}\cdots\lambda_p^{k_p/2}\dif\mu_n^{(p)}(\lambda_1,\ldots,\lambda_p).
\end{equation*}
By Theorem~\ref{thm:Bronk} and Lemma~\ref{lem:prodmarg}, it follows that
\begin{equation*}
\lim_{n\to\infty}\int\lambda_1^{k_1}\lambda_2^{k_2}\cdots\lambda_p^{k_p}\dif\sigma_n^{(p)}(\lambda_1,\ldots,\lambda_p)=
\prod_{i=1}^p\int\lambda_i^{k_i/2}\dif\nu_c(\lambda_i),
\end{equation*}
where $\nu_c$ is the free Poisson distribution of parameter $c$.
Therefore $\sigma_n^{(p)}$ converges in moments to $\btimes_1^p\rho$, where $\rho$ has density
\begin{equation*}
\frac{\dif\rho}{\dif t}=\frac{\sqrt{(d_1^2-t^2)(t^2-d_0^2)}}{\pi t}1_{[d_0,d_1]}(t),
\end{equation*}
with $d_0=1-\sqrt c$ and $d_1=1+\sqrt c$.
Now Theorem~\ref{thm:Rdiag} applies and finishes the proof.
\end{proof}

Every complex $n\times n$ matrix $A$ is unitarily conjugate to an upper triangular matrix: $A=USU^*$ where $U$ is unitary
and the $(i,j)$th entry of $S$ is zero if $i>j$.
If $A$ has $n$ distinct eigenvalues then the pair $(U,S)$ is unique up to replacement by $(UD,D^*SD)$, where $D$ is a diagonal unitary.
Given a random matrix $X\in\MEu_n$, one may ask for a corresponding random upper triangular matrix $S$ and random unitary matrix $U$
so that the distribution of $USU^*$ is equal to the distribution of $X$.
Then $X$ and $S$ will have the same $*$--moments with respect to the functional $\tau_n$.
For specificity, we may insist that the joint distribution of the pair $(U,S)$ be the same as the joint distribution of $(UD,D^*SD)$
for every non--random diagonal unitary $D$, (i.e.\ that the joint distribution of $(U,S)$ be invairant
under this action of the $n$--torus $\Tcirc^n$).
If the distribution of $X$ is invariant under conjugation by non--random unitaries and if $(U,S)$ is the pair of random matrices
as described above, then it is clear that the random unitary $U$ is distributed according to Haar measure on the $n\times n$ unitaries
and that $U$ and $S$ are independent.
In this case, the relavant question is only the distribution of $S$.
F.~Dyson answered this question when $X\in\GRM(n,\frac1n)$.
We state his result, and then make a slight modification to give, in conjunction with
Theorem~\ref{thm:RPmm}, an upper triangular matrix model for a circular free Poisson element.

\begin{thm}[Dyson, see~{\cite[A.35]{Mehta}}]
\label{thm:Dyson}
Let  $T(n)\in\UTGRM(n,\frac1n)$ and let $D(n)\in\MEu_n$ be a diagonal random matrix, whose diagonal
entries have joint density
\begin{equation}
\label{eq:UTdistr}
K^{(7)}_n\exp\Bigl(-n\sum_{i=1}^n|z_i|^2\Bigr)\prod_{1\le i<j\le n}|z_i-z_j|^2
\end{equation}
with respect to Lebesgue measure on $\Cpx^n$, for some constant $K^{(7)}_n$.
Let $U(n)\in\HURM(n)$,
suppose that $\big(D(n),T(n),U(n)\big)$ is an independent family of matrix--valued random variables and let
\begin{equation*}
X(n)=U(n)\big(D(n)+T(n)\big)U(n)^*.
\end{equation*}
Then $X(n)\in\GRM(n,\frac1n)$.
Consequently, $D(n)+T(n)$ converges in $*$--moments as $n\to\infty$ to a circular element.
\end{thm}

\begin{cor}
\label{cor:Dysondet}
Let $c\ge1$, let  $T(n)\in\UTGRM(n,\frac1n)$ and let $D_c(n)\in\MEu_n$ be a diagonal random matrix, whose diagonal
entries have joint density
\begin{equation*}
K^{(8)}_{c,n}\exp\Bigl(-n\sum_{i=1}^n|z_i|^2\Bigr)\Bigl(\prod_{i=1}^n|z_i|\Bigr)^{2(c-1)n}\prod_{1\le i<j\le n}|z_i-z_j|^2
\end{equation*}
with respect to Lebesgue measure on $\Cpx^n$, for some constant $K^{(8)}_{c,n}$.
Let $U(n)\in\HURM(n)$,
suppose that $\big(D_c(n),T(n),U(n)\big)$ is an independent family of matrix--valued random variables and let
\begin{equation}
\label{eq:Ynmat}
Y(n)=U(n)\big(D_c(n)+T(n)\big)U(n)^*.
\end{equation}
Then $Y(n)$ has density with respect to Lebesgue measure on $M_n(\Cpx)$ equal to
\begin{equation}
\label{eq:detdens}
K^{(1)}_{c,n}|\det Y|^{2(c-1)n}\exp\bigl(-n\Tr(Y^*Y)\bigr).
\end{equation}
Consequently, $D_c(n)+T(n)$ converges in $*$--moments as $n\to\infty$ to a circular free Poisson element of parameter $c$.
\end{cor}
\begin{proof}
Let $\Mc_n$ be the manifold of matrices in $M_n(\Cpx)$ having $n$ distinct eigenvalues.
Then $\Mc_n$ has full Lebesgue measure in $M_n(\Cpx)$.
Let $\UEu_n$ be the Lie group of $n\times n$ unitary matrices, and let $\TEu_n$ be the manifold of all upper triangular
$n\times n$ complex matrices, no two of whose diagonal elements are the same.
Let $\pi:\UEu_n\times\TEu_n\to\Mc_n$ be given by $\pi(U,S)=USU^*$.
Dyson proved his result by evaluating the Jacobian of $\pi$ (after throwing away the directions in $\ker\dif\pi$)
and thereby finding the measure $\sigma_n$ on $\TEu_n$ such that letting $\mu_n$ be Haar measure on $\UEu_n$, the push--forward
measure $\pi_*(\mu_n\times\sigma_n)$ on $\Mc_n$ has density $K_{1,n}^{(1)}\exp\big(-n\Tr(Y^*Y)\big)$ with respect to Lebesgue measure on $\Mc_n$,
i.e.\ the density of a random matrix $X(n)\in\GRM(n,\frac1n)$.
This measure $\sigma_n$ was found to have density
\begin{equation}
\label{eq:densS}
K_n^{(9)}\prod_{1\le i<j\le n}|S_{ii}-S_{jj}|^2\exp\big(-n\Tr(S^*S)\big)
\end{equation}
with respect to Lebesgue measure on $\TEu_n$, where for a matrix $S\in\TEu_n$, $S_{ii}$ is the $i$th diagonal entry of $S$;
this density~\eqref{eq:densS} is that of the matrix $D(n)+T(n)$ in Theorem~\ref{thm:Dyson}.

The matrix $D_c(n)+T(n)$ in the corollary has density
\begin{equation*}
K_{c,n}^{(10)}\prod_{1\le i<j\le n}|S_{ii}-S_{jj}|^2\exp\big(-n\Tr(S^*S)\big)|\det(S)|^{2(c-1)n}
\end{equation*}
with respect to Lebesgue measure on $\TEu_n$;
since $\det(USU^*)=\det(S)$, and building on Dyson's calculation, it follows that the random matrix $Y(n)$ of~\eqref{eq:Ynmat} has
density~\eqref{eq:detdens} with respect to Lebesgue measure on $\Mc_n$, as required.

An application of Theorem~\ref{thm:RPmm} shows that $Y(n)$, and hence also $D_c(n)+T(n)$, converges in $*$--moments as $n\to\infty$
to a circular free Poisson element.
\end{proof}

The following lemma shows that the diagonal entries of $D_c(n)$ are in a specific sense asymptotically independent.
This will allow their eventual decoupling; (see Remark~\ref{rem:*mom}).

\begin{lemma}
\label{lem:unifannu}
For $c\ge1$ and $n\in\Nats$ let $\mu_n$ be the probability measure on $\Cpx^n$ whose density with respect to Lebesgue measure is
\begin{equation*}
D_n(z_1,\ldots,z_n)=K^{(8)}_{c,n}\exp\Bigl(-n\sum_{i=1}^n|z_i|^2\Bigr)\Bigl(\prod_{i=1}^n|z_i|\Bigr)^{2(c-1)n}\prod_{1\le i<j\le n}|z_i-z_j|^2.
\end{equation*}
Given $p\in\{1,2,\ldots,n\}$ let $\mu_n^{(p)}$ be the marginal distribution of $\mu_n$ corresponding to the first $p$ variables $z_1,\ldots,z_p$.
Then for every $p\in\Nats$, $\mu_n^{(p)}$ converges in weak$^*$ topology and in $*$--moments as $n\to\infty$ to the product measure
$\btimes_1^p\rho$, where $\rho$ is uniform distribution on the annulus
$\{z\in\Cpx\mid\sqrt{c-1}<|z|<\sqrt c\}$.
\end{lemma}
\begin{proof}
This is quite similar to the proof of Lemma~\ref{lem:prodmarg}.
Let $\alpha=(c-1)n$.
Consider first the case $p=1$.
Let $\psi_0,\psi_1,\psi_2,\ldots$ be the polynomials obtained via Gram--Schmidt orthonormalization of the sequence $1,z,z^2,\ldots$
in $L^2\big(\Cpx,|z|^{2\alpha}e^{-n|z|^2}\dif(\RealPart z)\dif(\ImagPart z)\big)$.
Then
\begin{equation*}
\psi_k(z)=\sqrt{\frac{n^{k+\alpha+1}}{\pi\Gamma(k+\alpha+1)}}z^k.
\end{equation*}
Using the Vandermonde determinant we have
\begin{equation*}
\prod_{1\le i<j\le n}(z_j-z_i)=K_{c,n}^{(11)}\det\left(
\begin{matrix}
\psi_0(z_1)&\psi_0(z_2)&\cdots&\psi_0(z_n) \\[1ex]
\psi_1(z_1)&\psi_1(z_2)&\cdots&\psi_1(z_n) \\
\vdots&\vdots&&\vdots\\
\psi_{n-1}(z_1)&\psi_{n-1}(z_2)&\cdots&\psi_{n-1}(z_n)
\end{matrix}\right)
\end{equation*}
for some constant $K_{c,n}^{(11)}$.
Therefore
\begin{equation*}
D_n(z_1,\ldots,z_n)=K_{c,n}^{(12)}\left(\prod_{i=1}^n|z_i|\right)^{2\alpha}
\left|\sum_{\pi\in S_n}\sign(\pi)\prod_{i=1}^n\psi_{\pi(i)-1}(z_i)\right|^2
\exp\left(-n\sum_{i=1}^n|z_i|^2\right).
\end{equation*}
Writing
\begin{equation*}
\prod_{i=1}^n\psi_{\pi(i)-1}(z_i)=\psi_{\pi(1)-1}(z_1)\otimes\psi_{\pi(2)-1}(z_2)\otdt\psi_{\pi(n)-1}(z_n)
\end{equation*}
and noting that as $\pi$ ranges over the permutation group $S_n$ these form an orthonormal family
with respect to the measure
$\bigl(\prod_{i=1}^n|z_i|\bigr)^{2\alpha}\exp\bigl(-n\sum_{i=1}^n|z_i|^2\bigr)$ on $\Cpx^n$, we find $K_{c,n}^{(12)}=(n!)^{-1}$.
Moreover, the density of $\mu_n^{(1)}$ with respect to Lebesgue measure on $\Cpx$ is
\begin{align*}
D_{n,1}(z)&=\int_{\Cpx^{n-1}}D_n(z,z_2,\ldots,z_n)\dif(\RealPart z_1)\dif(\ImagPart z_1)\cdots\dif(\RealPart z_n)\dif(\ImagPart z_n) \\
&=\frac1n\left(\sum_{k=0}^{n-1}|\psi_k(z)|^2\right)|z|^{2\alpha}e^{-n|z|^2}
=\frac{n^\alpha}\pi\sum_{k=0}^{n-1}\frac{n^k|z|^{2k+2\alpha}}{\Gamma(k+\alpha+1)}e^{-n|z|^2}.
\end{align*}
We shall show that $\mu_n^{(1)}$ converges in $*$--moments to $\rho$.
Clearly if $a,b\in\Nats\cup\{0\}$ and if $a\ne b$ then
\begin{equation*}
\int_\Cpx z^a\zbar^{\,b}D_{n,1}(z)\dif(\RealPart z)\dif(\ImagPart z)
=0=\int_\Cpx z^a\zbar^{\,b}\dif\rho(z).
\end{equation*}
Hence we need only show
\begin{equation*}
\lim_{n\to\infty}\int_\Cpx|z|^{2b}\dif\mu_n^{(1)}(z)=\int_\Cpx|z|^{2b}\dif\rho(z)
\end{equation*}
for all $b\in\Nats\cup\{0\}$.
We have
\begin{align*}
\int_\Cpx|z|^{2b}\dif\mu_n^{(1)}(z)
&=\frac{n^{(c-1)n}}\pi\sum_{k=0}^{n-1}\frac{n^k}{\Gamma(k+(c-1)n+1)}\int_\Cpx|z|^{2(b+k+(c-1)n)}e^{-n|z|^2}\dif(\RealPart z)\dif(\ImagPart z) \\
&=\sum_{k=0}^{n-1}\frac{n^{k+(c-1)n}}{\Gamma(k+(c-1)n+1)}\int_0^\infty t^{b+k+(c-1)n}e^{-nt}\dif t.
\end{align*}
Writing
\begin{equation*}
f_n(t)=\sum_{k=0}^{n-1}\frac{n^{k+(c-1)n}}{\Gamma(k+(c-1)n+1)}t^{b+k+(c-1)n}e^{-nt},
\end{equation*}
we have
\begin{align*}
\int_\Cpx|z|^{2b}\dif\mu_n^{(1)}(z)&=\int_0^\infty t^bf_n(t)\dif t=-\frac1{b+1}\int_0^\infty t^{b+1}f_n'(t)\dif t= \displaybreak[0] \\[1ex]
&=\frac1{b+1}\left(\frac{\Gamma(cn+b+1)}{n^{b+1}\Gamma(cn)}-\frac{\Gamma((c-1)n+b+1)}{n^{b+1}\Gamma((c-1)n)}\right)= \displaybreak[1] \\[1ex]
&=\frac1{b+1}\Bigl(c(c+\tfrac1n)\cdots(c+\tfrac bn)-(c-1)\bigl((c-1)+\tfrac1n\bigr)\cdots\bigl((c-1)+\tfrac bn\bigr)\Bigr)  \\[1ex]
&\overset{n\to\infty}\longrightarrow\frac1{b+1}\bigl(c^{b+1}-(c-1)^{b+1}\bigr)=\int_\Cpx|z|^{2b}\dif\rho(z).
\end{align*}
Hence $\mu_n^{(1)}$ converges in $*$--moments to $\rho$ as $n\to\infty$;
since $\rho$ is compactly supported it follows that $\mu_n^{(1)}$ converges in the weak$^*$ topology to $\rho$.

The density of $\mu_n^{(2)}$ with respect to Lebesgue measure is
\begin{align*}
D_{n,2}(z_1,z_2)&=\int_{\Cpx^{n-2}}D_n(z_1,z_2,\ldots,z_n)\dif(\RealPart z_3)\dif(\ImagPart z_3)\cdots\dif(\RealPart z_n)\dif(\ImagPart z_n) \\
&=\begin{aligned}[t]
  \frac1{n!}\sum_{\pi,\sigma\in S_n}\sign(\pi)\sign(\sigma)
  \int_{\Cpx^{n-2}}\prod_{i=1}^n\Bigl(&\psi_{\pi(i)-1}(z_i)\overline{\psi_{\sigma(i)-1}(z_i)}|z_i|^{2\alpha}e^{-n|z_i|^2}\Bigr)\cdot \\
   &\cdot\dif(\RealPart z_3)\dif(\ImagPart z_3)\cdots\dif(\RealPart z_n)\dif(\ImagPart z_n) \end{aligned} \\
&=\frac1{n(n-1)}\sum_{\substack{0\le k,\ell\le n-1\\k\ne\ell}}
 \begin{aligned}[t]
  &\Bigl(|\psi_k(z_1)|^2|\psi_\ell(z_2)|^2-\psi_k(z_1)\psi_\ell(z_2)\overline{\psi_\ell(z_1)}\overline{\psi_k(z_2)}\Bigr)\cdot \\
  &\cdot|z_1|^{2\alpha}|z_2|^{2\alpha}e^{-n(|z_1|^2+|z_2|^2)} \end{aligned} \\
&=\frac n{n-1}D_{n,1}(z_1)D_{n,1}(z_2)-\frac1{n(n-1)}\left|\sum_{k=0}^{n-1}\psi_k(z_1)\overline{\psi_k(z_2)}\right|^2
 |z_1|^{2\alpha}|z_2|^{2\alpha}e^{-n(|z_1|^2+|z_2|^2)}.
\end{align*}
Hence as a linear functional on $C_0(\Cpx^2)$, the norm of $\mu_n^{(2)}-\mu_n^{(1)}\otimes\mu_n^{(1)}$ is bounded above by $2/(n-1)$.
Therefore $\mu_n^{(2)}$ converges in the weak$^*$  topology as $n\to\infty$ to $\rho\times\rho$.
Arguing as in the proof of Lemma~\ref{lem:prodmarg} and using Theorem~\ref{thm:ProdMeas}, we conclude that for every $p\ge1$,
$\mu_n^{(p)}$ converges in weak$^*$ topology to $\btimes_1^p\rho$, which we will denote by $\nu^{(p)}$.

It remains to show that $\mu_n^{(p)}$ converges to $\nu^{(p)}$ in $*$--moments, namely that
\begin{equation}
\label{eq:stm}
\lim_{n\to\infty}\int z_1^{k_1}\overline{z_1}^{\ell_1}\cdots z_p^{k_p}\overline{z_p}^{\ell_p}\dif\mu_n^{(p)}(z_1,\ldots,z_p)
=\int z_1^{k_1}\overline{z_1}^{\ell_1}\cdots z_p^{k_p}\overline{z_p}^{\ell_p}\dif\nu^{(p)}(z_1,\ldots,z_p)
\end{equation}
for every $k_1,\ldots,k_p,\ell_1,\ldots,\ell_p\in\Nats\cup\{0\}$.
Exactly as in the proof of Claim~\ref{claim:hpos}, one shows that if $h$ is a positive continuous function on $\Cpx^p$ then
\begin{equation*}
\liminf_{n\to\infty}\int h\dif\mu_n^{(p)}\ge\int h\dif\nu^{(p)}.
\end{equation*}
Then, considering the real and imaginary parts separately and arguing as in the proof of Claim~\ref{claim:mgfg}, one shows that
if $f$ and $g$ are continuous functions on $\Cpx^p$, if $g\ge0$, if $|f|\le g$ and if
$\lim_{n\to\infty}\int g\dif\mu_n^{(p)}=\int g\dif\nu^{(p)}$ then $\lim_{n\to\infty}\int f\dif\mu_n^{(p)}=\int f\dif\nu^{(p)}$.
But letting $m=k_1+\cdots+k_p+\ell_1+\cdots+\ell_p$, we have
\begin{equation*}
|z_1^{k_1}\overline{z_1}^{\ell_1}\cdots z_p^{k_p}\overline{z_p}^{\ell_p}|\le1+|z_1|^{2m}+\cdots+|z_p|^{2m}.
\end{equation*}
Moreover, because $\mu_n^{(1)}$ converges in $*$--moments to $\nu^{(1)}$, we have
\begin{align*}
\lim_{n\to\infty}\int(1+|z_1|^{2m}+\cdots+|z_p|^{2m})\dif\mu_n^{(p)}&=1+p\lim_{n\to\infty}\int|z_1|^{2m}\dif\mu_n^{(1)}
=1+p\int|z_1|^{2m}\dif\nu^{(1)}= \\
&=\int(1+|z_1|^{2m}+\cdots+|z_p|^{2m})\dif\nu^{(p)}.
\end{align*}
Hence we have~\eqref{eq:stm} and the lemma is proved.
\end{proof}

\begin{lemma}
\label{lemma:mupr}
Let $c\ge1$ and for every $n\in\Nats$ let $\mu_n$ and $\mu_n'$ be the probability measures on $\Cpx^n$
whose densities with respect to Lebesgue measure are, respectively,
\begin{align}
\label{eq:densmu}
D_n(z_1,\ldots,z_n)&=K_{c,n}^{(8)}\left(\prod_{i=1}^n|z_i|\right)^{2(c-1)n}\left(\prod_{1\le i<j\le n}|z_i-z_j|^2\right)
\exp\left(-n\sum_{i=1}^n|z_i|^2\right) \\
\label{eq:densmup}
D_n'(z_1,\ldots,z_n)&=K_{c,n}^{(13)}\left(\prod_{i=1}^n|z_i|\right)^{2(c-1)n}\left(\sum_{\pi\in S_n}\prod_{i=1}^n|z_i|^{2(\pi(i)-1)}\right)
\exp\left(-n\sum_{i=1}^n|z_i|^2\right).
\end{align}
For $p\in\{1,\ldots,n\}$ let $\mu_n^{(p)}$ and $(\mu_n')^{(p)}$ denote the marginal distributions of $\mu_n$ and, respectively, $\mu_n'$ corresponding
to the variables $z_1,\ldots,z_p$.
Then for every $p$, $(\mu_n')^{(p)}$ is obtained from $\mu_n^{(p)}$ by averaging over the action of the torus $\Tcirc^p$ on $\Cpx^p$
given by coordinate--wise multiplication:
\begin{equation*}
\Tcirc^p\times\Cpx^p\ni\bigl((w_1,\ldots,w_p),(z_1,\ldots,z_p)\bigr)\mapsto(w_1z_1,\ldots,w_pz_p).
\end{equation*}
Consequently, $(\mu_n')^{(p)}$ converges in $*$--moments and in weak$^*$ topology as $n\to\infty$ to the measure $\btimes_1^p\rho$,
where $\rho$ is the uniform distribution on the annulus $\{z\in\Cpx\mid\sqrt{c-1}<|z|<\sqrt c\}$.
\end{lemma}
\begin{proof}
Using the Vandermonde determinant we find
\begin{equation}
\label{eq:Vdm}
\prod_{1\le i<j\le n}|z_i-z_j|^2
=\sum_{\pi,\sigma\in S_n}\sign(\pi)\sign(\sigma)\prod_{i=1}^nz_i^{\pi(i)-1}\overline{z_i}^{\sigma(i)-1}.
\end{equation}
Averaging~\eqref{eq:Vdm} over the action of $\Tcirc^n$ gives
\begin{equation*}
\sum_{\pi\in S_n}\prod_{i=1}^n|z_i|^{2(\pi(i)-1)}.
\end{equation*}
From this one easily sees that$K_{c,n}^{(8)}=K_{c,n}^{(13)}$ and that $\mu_n'$ is obtained from $\mu_n$ by averaging.
In order to show that $(\mu_n')^{(p)}$ is obtained from $\mu_n^{(p)}$ by averaging, it suffices to note that for any
measure $\tau$ on $\Cpx^n$, the average over the action of $\Tcirc^p$ on the marginal distribution, $\tau^{(p)}$, corresponding to the first
$p$ variables, is equal to the marginal distribution of the average over the action of $\Tcirc^n$ on $\tau$.

Since, by Lemma~\ref{lem:unifannu}, $\mu_n^{(p)}$ converges in $*$--moments and in weak$^*$ topology as $n\to\infty$ to $\btimes_1^p\rho$, which
is invariant under the action of $\Tcirc^p$, it follows that $(\mu_n')^{(p)}$ converges in $*$--moments and in weak$^*$ topology to $\btimes_1^p\rho$.
\end{proof}

\begin{rem}
\label{rem:*mom}
Our main purpose in proving the immediately preceding two lemmas was to be able to conclude that
\begin{equation}
\label{eq:mommumup}
\lim_{n\to\infty}\left(\int z_1^{k_1}\overline{z_1}^{\ell_1}\cdots z_p^{k_p}\overline{z_p}^{\ell_p}\dif\mu_n^{(p)}
-\int z_1^{k_1}\overline{z_1}^{\ell_1}\cdots z_p^{k_p}\overline{z_p}^{\ell_p}\dif(\mu_n')^{(p)}\right)=0
\end{equation}
for every $p\in\Nats$ and every $k_1,\ldots,k_p,\ell_1,\ldots,\ell_p\in\Nats\cup\{0\}$.
It is possible to prove~\eqref{eq:mommumup} directly using the Vandermonde determinant and
combinatorial arguments, though this sort of proof is not
as satisfying as the one above involving Lemma~\ref{lem:unifannu}, where the limit measure is found.
\end{rem}

\begin{thm}
\label{thm:utfrPsn}
Let $c\ge1$ and $N\in\Nats$, and let $(A,\phi)$ be a W$^*$--noncommutative probability space with random variables $a_1,\ldots,a_N\in A$
and $b_{ij}\in A$ (\/$1\le i<j\le N$), where $a_j$ is a circular free Poisson element of parameter $(c-1)N+j$, where each $b_{ij}$ is a
circular element with $\phi(b_{ij}^*b_{ij})=1$, and where the family
\begin{equation*}
\bigl((\{a_j^*,a_j\})_{1\le j\le N},(\{b_{ij}^*,b_{ij}\})_{1\le i<j\le N}\bigr)
\end{equation*}
is free.
Consider the W$^*$--noncommutative probability space $(M_N(A),\phi_N)$, where
\begin{equation*}
\phi_N\bigl((x_{ij})_{1\le i,j\le N}\bigr)=N^{-1}\sum_{j=1}^N\phi(x_{jj}),
\end{equation*}
and consider the random variable
\begin{equation*}
x=\frac1{\sqrt N}\left(\begin{matrix}a_1&b_{12}&b_{13}&\cdots&b_{1,N-1}&b_{1N}\\0&a_2&b_{23}&\cdots&b_{2,N-1}&b_{2N}\\0&0&a_3&\ddots&\vdots&b_{3N}\\
\vdots&\vdots&\ddots&\ddots&\ddots&\vdots\\ 0&0&\cdots&0&a_{N-1}&b_{N-1,N}\\0&0&\cdots&0&0&a_N\end{matrix}\right)\in M_N(A).
\end{equation*}
Then $x$ is a circular free Poisson element of parameter $c$.
\end{thm}
\begin{proof}
For every $n\in\Nats$ let $Y(nN)$ be an $nN\times nN$ random matrix whose distribution has density with respect to Lebesgue measure
\begin{equation*}
K_{c,nN}^{(1)}|\det(Y)|^{2(c-1)nN}\exp\bigl(-nN\Tr(Y^*Y)\bigr).
\end{equation*}
Then by Theorem~\ref{thm:RPmm}, $Y(nN)$ converges in $*$--moments as $n\to\infty$ to a circular free Poisson element of parameter $c$.
By Corollary~\ref{cor:Dysondet}, each $Y(nN)$ has the same $*$--moments as $S^{(1)}(nN)\eqdef D^{(1)}(nN)+T(nN)$,
where $T(nN)\in\UTGRM(nN,\frac1{nN})$, $D^{(1)}(nN)$
is a diagonal $nN\times nN$ random matrix, the distribution of whose diagonal entries has density
\begin{equation}
\label{eq:k=1}
K^{(8)}_{c,nN}\left(\prod_{i=1}^{nN}|z_i|\right)^{2(c-1)nN}\left(\prod_{1\le i<j\le nN}|z_i-z_j|^2\right)\exp\left(-nN\sum_{i=1}^n|z_i|^2\right)
\end{equation}
with respect to Lebesgue measure on $\Cpx^n$, and where $D^{(1)}(nN)$ and $T(nN)$ are independent.
We will use previous results to show that also each $S^{(k)}(nN)\eqdef D^{(k)}(nN)+T(nN)$, ($k\in\{2,3,4,5,6\}$),
converges in  $*$--moments as $n\to\infty$ to a circular free Poisson element of parameter $c$, where $D^{(k)}(nN)$ is a diagonal random matrix
such that $D^{(k)}(nN)$ and $T(nN)$ are independent and where the joint distributions of the diagonal entries of $D^{(k)}(nN)$ have the
following densities with respect to Lebesgue measure on $\Cpx^n$:
\begin{align*}
\text{for }k=2:
\quad&K_{c,nN}^{(13)}\left(\prod_{i=1}^{nN}|z_i|\right)^{2(c-1)nN}\left(\sum_{\pi\in S_{nN}}\prod_{i=1}^{nN}|z_i|^{2(\pi(i)-1)}\right)
\exp\left(-nN\sum_{i=1}^{nN}|z_i|^2\right) \displaybreak[0] \\
\text{for }k=3:
\quad&K_{c,nN}^{(14)}\left(\prod_{i=1}^{nN}|z_i|\right)^{2(c-1)nN}\left(\prod_{i=1}^{nN}|z_i|^{2(i-1)}\right)
\exp\left(-nN\sum_{i=1}^{nN}|z_i|^2\right) \displaybreak[0] \\
\text{for }k=4:
\quad&
\begin{aligned}[t]K_{c,nN}^{(15)}\prod_{j=1}^N\Biggl(\left(\prod_{i=1}^n|z_{(j-1)n+i}|\right)^{2((c-1)N+j-1)n}&\left(\prod_{i=1}^n|z_{(j-1)n+i}|^{2(i-1)}\right)\cdot \\
&\quad\cdot\exp\left(-nN\sum_{i=1}^n|z_{(j-1)n+i}|^2\right)\Biggr)\end{aligned} \displaybreak[0] \\
\text{for }k=5:
\quad&
\begin{aligned}[t]K_{c,nN}^{(16)}\prod_{j=1}^N\Biggl(\left(\prod_{i=1}^n|z_{(j-1)n+i}|\right)^{2((c-1)N+j-1)n}
&\left(\sum_{\pi\in S_n}\prod_{i=1}^n|z_{(j-1)n+i}|^{2(\pi(i)-1)}\right)\cdot \\
&\qquad\cdot\exp\left(-nN\sum_{i=1}^n|z_{(j-1)n+i}|^2\right)\Biggr)\end{aligned} \displaybreak[0] \\
\text{for }k=6:
\quad&
\begin{aligned}[t]K_{c,nN}^{(17)}\prod_{j=1}^N\Biggl(\left(\prod_{i=1}^n|z_{(j-1)n+i}|\right)^{2((c-1)N+j-1)n}
&\left(\prod_{1\le i<i'\le n}|z_{(j-1)n+i}-z_{(j-1)n+i'}|^2\right)\cdot \\
&\qquad\cdot\exp\left(-nN\sum_{i=1}^n|z_{(j-1)n+i}|^2\right)\Biggr).\end{aligned}
\end{align*}
The proof that $S^{(k)}(nN)$ converges in $*$--distribution to a circular free Poisson element relies for $k=2$ on
Lemmas~\ref{lem:unifannu} and~\ref{lemma:mupr}, (see Remark~\ref{rem:*mom}), and Theorem~\ref{thm:utgrm};
for $k=3$ we use Theorem~\ref{thm:*momPerm};
the density for $k=4$ is just a rewriting of that for $k=3$;
for $k=5$ we use again Theorem~\ref{thm:*momPerm};
for $k=6$ we use again Lemmas~\ref{lem:unifannu} and~\ref{lemma:mupr}, and Theorem~\ref{thm:utgrm}.

We may characterize the above successive transformations as follows:
from~\eqref{eq:k=1} to $k=2$ is decoupling;
from $k=2$ to $k=3$ is desymmetrization;
from $k=3$ to $k=4$ is regrouping;
from $k=4$ to $k=5$ is partial resymmetrization;
from $k=5$ to $k=6$ is partial recoupling.

Taking blocks of consequetive rows and columns to write  $D^{(6)}(nN)+T(nN)$ as an $N\times N$ matrix of $n\times n$ random matrices, we have
\begin{equation*}
S^{(6)}(nN)=\frac1{\sqrt N}\left(
\begin{matrix}
A^{(6)}_1&B^{(6)}_{12}&B^{(6)}_{13}&\cdots&B^{(6)}_{1,N-1}&B^{(6)}_{1N}\\
0&A^{(6)}_2&B^{(6)}_{23}&\cdots&B^{(6)}_{2,N-1}&B^{(6)}_{2N}\\
0&0&A^{(6)}_3&\ddots&\vdots&B^{(6)}_{3N}\\
\vdots&\vdots&\ddots&\ddots&\ddots&\vdots\\
0&0&\cdots&0&A^{(6)}_{N-1}&B^{(6)}_{N-1,N}\\
0&0&\cdots&0&0&A^{(6)}_N\end{matrix}
\right),
\end{equation*}
where
\begin{equation*}
\Bigl(\bigl(A_j^{(6)}(n)\bigr)_{1\le j\le N},\bigl(B_{ij}^{(6)}(n)\bigr)_{1\le i<j\le N}\Bigr)
\end{equation*}
is an independent family of matrix--valued random variables,
where $B_{ij}^{(6)}(n)\in\GRM(n,\frac1n)$ for every $1\le i<j\le N$ and where $A_j^{(6)}(n)=D_j^{(6)}(n)+T_j(n)$
with $T_j(n)\in\UTGRM(n,\frac1n)$, with $D_j^{(6)}(n)$ a diagonal random matrix, the joint distribution of whose diagonal entries
has density with respect to Lebesgue measure
\begin{equation*}
K_{(c-1)N+j,n}^{(8)}\left(\prod_{i=1}^n|z_i|\right)^{2((c-1)N+j-1)n}\left(\prod_{1\le i<i'\le n}^n|z_i-z_{i'}|^2\right)\exp\left(-n\sum_{i=1}^n|z_i|^2\right)
\end{equation*}
and with $D_j^{(6)}(n)$ and $T_j(n)$ independent.

Let $U_j(n)\in\HURM(n)$, ($1\le j\le N$), be such that 
\begin{equation*}
\Bigl(\bigl(A_j^{(6)}(n)\bigr)_{1\le j\le N},\bigl(B_{ij}^{(6)}(n)\bigr)_{1\le i<j\le N},\bigl(U_j(n)\bigr)_{1\le j\le N}\Bigr)
\end{equation*}
is an independent family of matrix--valued random variables.
By conjugating the matrix $S^{(6)}(nN)$ with $\diag(U_1(n),U_2(n),\ldots,U_N(n))$ and by
using Corollary~\ref{cor:Dysondet} and the fact that the class $\GRM(n,1/n)$ is invariant
under left and right multiplication by independent unitaries, it follows that
$S^{(6)}(nN)$ has the same $*$--moments, as
\begin{equation*}
S^{(7)}(nN)=\frac1{\sqrt N}\left(
\begin{matrix}
A^{(7)}_1&B^{(7)}_{12}&B^{(7)}_{13}&\cdots&B^{(7)}_{1,N-1}&B^{(7)}_{1N}\\
0&A^{(7)}_2&B^{(7)}_{23}&\cdots&B^{(7)}_{2,N-1}&B^{(7)}_{2N}\\
0&0&A^{(7)}_3&\ddots&\vdots&B^{(7)}_{3N}\\
\vdots&\vdots&\ddots&\ddots&\ddots&\vdots\\
0&0&\cdots&0&A^{(7)}_{N-1}&B^{(7)}_{N-1,N}\\
0&0&\cdots&0&0&A^{(7)}_N\end{matrix}
\right),
\end{equation*}
where
\begin{equation}
\label{eq:AB7}
\Bigl(\bigl(A_j^{(7)}(n)\bigr)_{1\le j\le N},\bigl(B_{ij}^{(7)}(n)\bigr)_{1\le i<j\le N}\Bigr)
\end{equation}
is an independent family of matrix--valued random variables,
where $B_{ij}^{(7)}(n)\in\GRM(n,\frac1n)$ for every $1\le i<j\le N$ and where the distribution of $A_j^{(7)}(n)$ has density
\begin{equation*}
K_{(c-1)N+j,n}^{(1)}|\det(A)|^{2((c-1)N+j-1)n}\exp(-n\Tr(A^*A))
\end{equation*}
with respect to Lebesgue measure on $M_n(\Cpx)$.

If $(V_j^{(1)})_{1\le j\le N}$, $(V_j^{(2)})_{1\le j\le N}$, $(U_{ij}^{(1)})_{1\le i<j\le N}$, $(U_{ij}^{(2)})_{1\le i<j\le N}$, are non--random
$n\times n$ unitary matrices, then
\begin{equation}
\label{eq:UVAB7}
\Bigl(\bigl(V_j^{(1)}A_j^{(7)}(n)V_j^{(2)}\bigr)_{1\le j\le N},\bigl(U_{ij}^{(1)}B_{ij}^{(7)}(n)U_{ij}^{(2)}\bigr)_{1\le i<j\le N}\Bigr)
\end{equation}
continues to be an independent family of matrix--valued random variables, $V_j^{(1)}A_j^{(7)}(n)V_j^{(2)}$ has the same distribution
as $A_j^{(7)}$ and $U_{ij}^{(1)}B_{ij}^{(7)}(n)U_{ij}^{(2)}$ has the same distribution as $B_{ij}^{(7)}(n)$.
Therefore, the family~\eqref{eq:UVAB7} has the same joint $*$--moments as the family~\eqref{eq:AB7}.
Taking into account also Theorem~\ref{thm:Bronk} and Lemma~\ref{lem:prodmarg} (as in the proof of Theorem~\ref{thm:RPmm}),
we see that the conditions of Theorem~\ref{thm:Rdiag} are fulfilled, allowing us to conclude that the family~\eqref{eq:AB7} is asymptotically $*$--free
as $n\to\infty$.
Moreover, (by Theorem~\ref{thm:RPmm}), each $A_j^{(7)}(n)$ converges in $*$--moments to a circular free Poisson element of parameter $(c-1)N+j$,
while $B_{ij}^{(7)}(n)$ converges in $*$--moments to a circular element.
Therefore, the entries of the matrix $S^{(7)}(n)$ model as $n\to\infty$ the entries of the matrix $x$ in the statement of the theorem.
As $S^{(7)}(n)$ converges in $*$--moments to a circular free Poisson element of parameter $c$, the theorem is proved.
\end{proof}

\section{Invariant subspaces for a circular free Poisson element}
\label{sec:iscfP}

In this section,
we will apply Theorem~\ref{thm:utfrPsn} and the general results of~\S\ref{sec:isut} to exhibit invariant subspaces for a circular free Poisson element.
We will rely on the result of Haagerup and Larsen~\cite[Example 5.2]{HL} that the spectrum of a circular free Poisson element of parameter $c$ is
$\{z\in\Cpx\mid\sqrt{c-1}\le|z|\le\sqrt c\}$.

\begin{thm}
\label{thm:fPinv}
Let $(\Mc,\psi)$ be a W$^*$--noncommutative probability space with $\psi$ faithful, let $c\ge1$ and let $y\in\Mc$ be a circular free Poisson element
of parameter $c$.
Given $r\ge0$, let $p_r(y)\in\Mc$ be the projection onto the invariant subspace of $y$ as in Definition~\ref{lemmadef:invvN}.
Then
\begin{equation}
\label{eq:psipry}
\psi(p_r(y))=\begin{cases}
0&\text{if }r\le\sqrt{c-1}\\
r^2-(c-1)&\text{if }\sqrt{c-1}\le r\le\sqrt c \\
1&\text{if }r\ge\sqrt c.
\end{cases}
\end{equation}
\end{thm}
\begin{proof}
We may without loss of generality assume that $\Mc=\{y\}''$, which implies $\Mc\cong L(F_2)$ and $\psi$ is a trace.
Let $N\in\Nats$ and let
\begin{equation*}
x=\frac1{\sqrt N}\left(\begin{matrix}a_1&b_{12}&b_{13}&\cdots&b_{1,N-1}&b_{1N}\\0&a_2&b_{23}&\cdots&b_{2,N-1}&b_{2N}\\0&0&a_3&\ddots&\vdots&b_{3N}\\
\vdots&\vdots&\ddots&\ddots&\ddots&\vdots\\ 0&0&\cdots&0&a_{N-1}&b_{N-1,N}\\0&0&\cdots&0&0&a_N\end{matrix}\right)\in M_N(A).
\end{equation*}
be the circular free Poisson element of parameter $c$ as in Theorem~\ref{thm:utfrPsn}, where we take $(A,\phi)$ to be a W$^*$--noncommutative probability space.
Thus $a_j$ is a circular free Poisson element of parameter $(c-1)N+j$, each $b_{ij}$ is a circular element and the collection of all $a_j$ and $b_{ij}$ is $*$--free.
For $k\in\{1,\ldots,N\}$, let
\begin{equation*}
e_k=\diag(\,\underset{k\text{ times}}{\underbrace{1,\ldots,1}},0,\ldots,0)\in M_n(\Cpx1)\subseteq M_n(A).
\end{equation*}
Another application of Theorem~\ref{thm:utfrPsn} shows that in the W$^*$--noncommutative probability space
\begin{equation*}
\bigl(e_kM_N(A)e_k,\sqrt{\tfrac Nk}\phi_N\restrict_{e_kM_N(A)e_k}\bigr)\cong(M_k(A),\phi_k),
\end{equation*}
the element $\sqrt{\frac Nk}e_kxe_k$ is a circular free Poisson element of parameter $\frac Nk(c-1)+1$.
Hence by~\cite[Example 5.2]{HL}, $e_kxe_k$ has spectrum
\begin{equation*}
\Bigl\{z\in\Cpx\Bigm|\sqrt{c-1}\le|z|\le\sqrt{(c-1)+\tfrac kN}\,\Bigr\}.
\end{equation*}
Similarly, if $k<N$ then denoting by $1_N$ the identity element of $M_N(A)$, we find that in the W$^*$--noncommutative probability space
\begin{equation*}
\bigl((1_N-e_k)M_N(A)(1_N-e_k),\sqrt{\tfrac N{N-k}}\Phi_N\restrict_{(1_N-e_k)M_N(A)(1_N-e_k)}\bigr)\cong(M_{N-k}(A),\phi_{N-k})
\end{equation*}
the element $\sqrt{\frac N{N-k}}(1_N-e_k)x(1_N-e_k)$ is a circular free Poisson element of parameter $\frac N{N-k}c$.
Hence $(1_N-e_k)x(1_N-e_k)$ has spectrum
\begin{equation*}
\Bigl\{z\in\Cpx\Bigm|\sqrt{(c-1)+\tfrac kN}\le|z|\le\sqrt c\,\Bigr\}.
\end{equation*}
Therefore, by Proposition~\ref{prop-uptrinv}, if $k\le N-2$ and if
\begin{equation*}
\sqrt{(c-1)+\tfrac kN}\le r<\sqrt{(c-1)+\tfrac{k+1}N}
\end{equation*}
then $e_k\le p_r(x)\le e_{k+1}$ and consequently $\frac kN\le\psi(p_r(y))\le\frac{k+1}N$.
Letting $N$ grow without bound and choosing $k$ appropriately implies~\eqref{eq:psipry}.
\end{proof}

Some further facts concerning these projections $p_r(y)$ are collected below in Theorem~\ref{thm:ypy}, for the proof of which we will use
the following lemma.
\begin{lemma}
\label{lemma:yconv}
For every $c\ge1$ let $(A_c,\phi_c)$ be a W$^*$--noncommutative probability space and let $y_c\in A$ be a circular free Poisson element of parameter $c$.
Given $c_0\ge1$ and a sequence $(c_n)_1^\infty$ in $[1,\infty)$ converging to $c_0$, we have that
$y_{c_n}$ converges in $*$--moments to $y_{c_0}$ as $n\to\infty$.
\end{lemma}
\begin{proof}
The positive part $h_c$
of $y_c$ has the same moments as the measure $\rho_c$ on $\Reals_+$ whose density with respect to Lebesgue measure is
\begin{equation*}
\frac{\dif\rho_c}{\dif t}=\frac{\sqrt{(d_1^2-t^2)(t^2-d_0^2)}}{\pi t}1_{[d_0,d_1]}(t),
\end{equation*}
with $d_0=1-\sqrt c$ and $d_1=1+\sqrt c$.
Since $y_c$ has the polar decomposition $y_c=u_ch_c$ where $u_c$ is a Haar unitary and where $u_c$ and $h_c$ are $*$--free,
the $*$--moments of $y_c$ can be expressed as certain polynomials in the moments of $\rho_c$.
Clearly, the $k$th moment of $\rho_{c_n}$ converges to the $k$th moment of $\rho_{c_0}$ as $n\to\infty$.
\end{proof}

\begin{thm}
\label{thm:ypy}
Let $y$ be a circular free Poisson element of parameter $c$ in some W$^*$--noncommutative probability space $(\Mc,\psi)$, with $\psi$ faithful.
Then
\renewcommand{\labelenumi}{(\roman{enumi})}
\begin{enumerate}
\item $p_s(y)$ converges in the strong$^*$ topology to $p_r(y)$ as $s\to r$.
\vskip1ex
\item If $\sqrt{c-1}<r<\sqrt c$ then in the W$^*$--noncommutative probability space
\begin{equation}
\label{eq:ncpspr}
\bigg(p_r(y)\Mc p_r(y),\tfrac1{\psi(p_r(y))}\psi\restrict_{p_r(y)\Mc p_r(y)}\bigg),
\end{equation}
$\psi(p_r(y))^{-1/2}yp_r(y)$
is a circular free Poisson element of parameter $1+(c-1)/\psi(p_r(y))$.
Hence the spectrum of $yp_r(y)$ relative to $p_r(y)\Mc p_r(y)$ is
\begin{equation*}
\sigma(yp_r(y))=\{z\in\Cpx\mid\sqrt{c-1}\le|z|\le\sqrt r\}.
\end{equation*}
\item If $\sqrt{c-1}<r<\sqrt c$ then in the W$^*$--noncommutative probability space
\begin{equation*}
\bigg((1-p_r(y))\Mc(1-p_r(y)),\tfrac1{1-\psi(p_r(y))}\psi\restrict_{(1-p_r(y))\Mc(1-p_r(y))}\bigg),
\end{equation*}
$\big(1-\psi(p_r(y))\big)^{-1/2}(1-p_r(y))y$ is a circular free Poisson element of parameter $c/\big(1-\psi(p_r(y))\big)$.
Hence the spectrum of $(1-p_r(y))y$ relative to $(1-p_r(y))\Mc(1-p_r(y))$ is
\begin{equation*}
\sigma((1-p_r(y))y)=\{z\in\Cpx\mid\sqrt r\le|z|\le\sqrt c\}.
\end{equation*}
\item If $x$ is the upper triangular $N\times N$ matrix given in Theorem~\ref{thm:utfrPsn} that is a circular free Poisson element of parameter $c$,
then for every $k\in\{0,1,,\ldots,N\}$ and letting $r=\sqrt{(c-1)+k/N}$, we have $p_r(x)=\diag(\,\underset{k\text{ times}}{\underbrace{1,\ldots,1}},0,\ldots,0)$.
\end{enumerate}
\end{thm}
\begin{proof}
We know from general principles that $p_{r'}(y)\le p_r(y)$ if $r'<r$, and from Theorem~\ref{thm:fPinv} we have that $\lim_{s\to r}\psi(p_s(y))=\psi(p_r(y))$;
as $\psi$ is faithful we conclude~(i).

Let us now prove~(iv).
Arguing as in the proof of Theorem~\ref{thm:fPinv}, we have $e_k\le p_r(x)$ whenever $r>\sqrt{(c-1)+k/N}$.
We may take the W$^*$-noncommutative probability space $(A,\phi)$ so that $\phi$ is a faithful trace,
in which case, since $\inf\{\psi(p_r(x))\mid r>\sqrt{(c-1)+k/N}\}=k/N$, it follows that
\begin{equation*}
e_k=\bigwedge\Big\{p_r(x)\,\Big|\, r>\sqrt{(c-1)+k/N}\Big\}.
\end{equation*}
Thus $e_k$ is the limit in strong$^*$ topology of $p_r(x)$ as $r$ tends to $\sqrt{(c-1)+k/N}$ from above.
Using~(i), it follows that $e_k=p_{\sqrt{(c-1)+k/N}}(x)$.

For~(ii), let us show that $\psi(p_r(y))^{-1/2}yp_r(y)$ is circular free Poisson of the desired parameter, first in the case when $\psi(p_r(y))=k/N$ is rational.
We may take $(\Mc,\psi)$ to be $(M_N(A),\phi_N)$ and $y$ to be equal to the $N\times N$ matrix $x$ as in Theorem~\ref{thm:utfrPsn}.
By~(iv), the noncommutative probability space~\eqref{eq:ncpspr} is $(M_k(A),\phi_k)$ and
\begin{equation*}
\frac1{\sqrt{\psi(p_r(y))}}yp_r(y)=\frac1{\sqrt k}\left(\begin{matrix}a_1&b_{12}&b_{13}&\cdots&b_{1,k-1}&b_{1k}\\0&a_2&b_{23}&\cdots&b_{2,k-1}&b_{2k}\\0&0&a_3&\ddots&\vdots&b_{3k}\\
\vdots&\vdots&\ddots&\ddots&\ddots&\vdots\\ 0&0&\cdots&0&a_{k-1}&b_{k-1,k}\\0&0&\cdots&0&0&a_k\end{matrix}\right),
\end{equation*}
where $a_j$ is circular free Poisson of parameter $(c-1)N+j$.
Applying again Theorem~\ref{thm:utfrPsn}, we obtain that $\psi(p_r(y))^{-1/2}yp_r(y)$ is circular free Poisson of parameter $1+(c-1)/\psi(p_r(y))$.
When $r$ is such that $\psi(p_r(y))$ is irrational,
then using~(i) we have that $yp_r(y)$ is the strong$^*$ limit of $yp_s(y)$ as $s$ tends to $r$ through rational numbers.
Hence by Lemma~\ref{lemma:yconv} and the continuity in $r$ of of $\psi(p_r(y))$ implied by Theorem~\ref{thm:utfrPsn},
it follows that $\psi(p_r(y))^{-1/2}yp_r(y)$ is circular free Poisson of parameter $1+(c-1)/\psi(p_r(y))$.
The statement about the spectrum follows from the result of Haagerup and Larsen~\cite{HL} that we've been using repeatedly.

Part~(iii) is proved similarly.
When $\psi(p_r(y))=k/N$ is rational then we get
\begin{equation*}
\frac1{\sqrt{1-\psi(p_r(y))}}(1-p_r(y))y=
\frac1{\sqrt{N-k}}\left(
\begin{matrix}a_{k+1}&b_{k+1,k+2}&b_{k+1,k+3}&\cdots&b_{k+1,N-1}&b_{k+1,N}\\
0&a_{k+2}&b_{k+2,k+3}&\cdots&b_{k+2,N-1}&b_{k+2,N}\\
0&0&a_{k+3}&\ddots&\vdots&b_{k+3,N}\\
\vdots&\vdots&\ddots&\ddots&\ddots&\vdots\\
0&0&\cdots&0&a_{N-1}&b_{N-1,N}\\0&0&\cdots&0&0&a_N\end{matrix}\right),
\end{equation*}
where $a_{k+j}$ is circular free Poisson of parameter $(c-1)N+k+j$.
The remaining part of the argument is like for~(ii) above.
\end{proof}

The following proposition shows that $p_r(y)$ is characterized by the spectral conditions in~(ii) and~(iii) of Theorem~\ref{thm:ypy}.
\begin{prop}
\label{prop:specchar}
Let $Y$ be a circular free Poisson element of parameter $c\ge1$ in a W$^*$--probability space $(\Mc,\psi)$,
with $\psi$ faithful, and let $\sqrt{c-1}<r<\sqrt c$.
Suppose $p\in\Mc$ is a projection such that
\renewcommand{\labelenumi}{(\roman{enumi})}
\begin{enumerate}
\item $yp=pyp$
\item $\sigma_{p\Mc p}(yp)\subseteq\{z\in\Cpx\mid\sqrt{c-1}\le|z|\le r\}$
\item $\sigma_{(1-p)\Mc(1-p)}((1-p)y)\subseteq\{z\in\Cpx\mid r\le|z|\le\sqrt c\}$.
\end{enumerate}
Then $p=p_r(y)$.
\end{prop}
\begin{proof}
Note that~(i) implies $(1-p)y=(1-p)y(1-p)$.
Let $\Mc$ be normally and faithfully represented on a Hilbert space $\HEu$.
For $\xi\in p\HEu$ we have
\begin{equation*}
\limsup_{n\to\infty}\|y^n\xi\|^{1/n}\le\limsup_{n\to\infty}\|(pyp)^n\|^{1/n}\le r,
\end{equation*}
where the last inequality is because the spectral radius of $pyp$ is $\le r$.
Hence $p\le p_r(y)$.

In order to prove the reverse inequality, it will suffice to show $p\ge p_s(y)$ for all $0\le s<r$,
because $s\mapsto p_s(y)$ is strong$^*$--continuous by Theorem~\ref{thm:ypy}(i).
Let $0\le s<r$, $\xi\in(1-p)\HEu$ and let $\eta\in E_s(y)$, i.e.\ 
\begin{equation}
\label{eq:yneta}
\limsup_{n\to\infty}\|y^n\eta\|^{1/n}\le s.
\end{equation}
Set $\xi_n=\big((1-p)y^*(1-p)\big)^{-n}\xi$.
Then $\xi_n\in(1-p)\HEu$ and, because the spectral radius of $\big((1-p)y^*(1-p)\big)^{-1}$ is $\le1/r$, we have
$\limsup_{n\to\infty}\|\xi_n\|^{1/n}\le1/r$.
Since $(1-p)\HEu$ is an invariant subspace for $y^*$, we have $\xi=\big((1-p)y^*(1-p)\big)^n\xi_n=(y^*)^n\xi_n$.
Therefore $\langle\xi,\eta\rangle=\langle(y^*)^n\xi_n,\eta\rangle=\langle\xi_n,y^n\eta\rangle$,
so using~\eqref{eq:yneta} and Schwarz's inequality we have $\limsup_{n\to\infty}|\langle\xi,\eta\rangle|^{1/n}\le s/r<1$,
which shows that $\langle\xi,\eta\rangle=0$.
Hence $(1-p)\HEu\perp\overline{E_s(y)}=p_s(y)\HEu$ and therefore $p_s(y)\le p$.
\end{proof}

\begin{rem}
The proof above shows that the subspace
\begin{equation*}
E_r(y)=\{\xi\in\HEu\mid\limsup_{k\to\infty}\nm{y^k\xi}^{1/k}\le r\}
\end{equation*}
is closed; thus we have $p_r(y)\HEu=E_r(y)$, without taking the closure.
\end{rem}

The next example, however, shows that the sort of spectral decomposition found in Theorem~\ref{thm:ypy} and closedness of the subspace $E_r(y)$
do not always hold.

\begin{example}
{\rm Let $\HEu=\bigoplus_{k=2}^\infty\HEu_k$, where $\HEu_k$ is $k$--dimensional Hilbert space with orthonormal basis $e^{(k)}_1,\ldots,e^{(k)}_k$
and let $T=\bigoplus_{k=2}^\infty T_k$, where $T_k\in B(\HEu_k)$ is the nilpotent operator
\begin{equation*}
T_ke^{(k)}_j=\cases 0 & \quad j=1 \\ e_{j-1}^{(k)} & \quad2\le j\le k. \endcases
\end{equation*}
It is well-known that the spectrum of $T$ is the closed unit disk
$\overline{\Db}$ --- see for example Brown~\cite[Example~4.10]{Brown}. 
However, if $r>0$ then $E_r(T)$ is dense in $\HEu$, so $p_r(T)=1$ and the spectrum of $Tp_r(T)$ is $\overline{\Db}$.
Moreover, if $0<r<1$ then the vector $\sum_{k=2}^\infty\frac1k e^{(k)}_k$ is not an element of $E_r(T)$;
this shows that $E_r(T)$ is not closed.
Note that $T\in\bigoplus_{k=2}^{\infty}B(\HEu_k)$, which is a finite von
Neumann algebra.
}
\end{example}

\end{spacing}
\begin{spacing}{1.0}

\vskip1ex
\scriptsize
\noindent{\sc Ken Dykema \\ Department of Mathematics, Texas A\&M University, College Station TX 77843--3368, USA} \\
\noindent{\sl E-mail:} {\tt Ken.Dykema@math.tamu.edu} \\
\noindent{\sl Internet URL:} {\tt http://www.math.tamu.edu/$\,\widetilde{\;}$Ken.Dykema/} \\
\vskip1ex
\noindent{\sc Uffe Haagerup* \\ Department of Mathematics and Computer Science, \\ University of Southern Denmark --- Odense University, Campusvej 55, 5230 Odense M, Denmark} \\
\noindent{\sl E-mail:} {\tt haagerup@imada.sdu.dk} \\
\noindent{\sl Internet URL:} {\tt http://www.imada.sdu.dk/$\,\widetilde{\;}$haagerup/} \\

\vskip1ex

\noindent *MaPhySto, Centre for Mathematical Physics and
Stochastics, funded by a grant from The Danish National Research Foundation.

\end{spacing}
\end{document}